%% file: main.tex
\documentclass[11pt]{article}
\input{preamble}

\title{On the Chv\'atal--Gomory Closure of Integer Quadratic Programs}
\author{%
   Linchuan Wei\thanks{%
    Corresponding author. School of Data Science,
    The Chinese University of Hong Kong, Shenzhen,
    2001 Longxiang Road, Longgang District,
    Shenzhen, Guangdong 518172, China.
    Email: \href{mailto:linchuanwei@cuhk.edu.cn}{linchuanwei@cuhk.edu.cn}.
    ORCID: \href{https://orcid.org/0000-0003-1522-4501}{0000-0003-1522-4501}.}
  \and
  Tong Xu\thanks{%
    Industrial Engineering \& Management Science, Northwestern University, 2145 Sheridan Road, Evanston, Illinois 60208, U.S.A.
    Email: tongxu2027@u.northwestern.edu.
    ORCID: \href{https://orcid.org/0009-0008-2654-5152}{0009-0008-2654-5152}.}
  \and
  Kim Yu \thanks{%
    D\'epartement d'informatique et de recherche op\'erationnelle,
    Universit\'e de Montr\'eal,
    2920 Tour Rd, Montr\'eal, Quebec H3T 1N8, Canada.
    Email: kimyu@iro.umontreal.ca.
    ORCID: \href{https://orcid.org/0000-0002-0609-3645}{0000-0002-0609-3645}.}
}
\date{}

\begin{document}
\maketitle

\begin{abstract}
The integer quadratic moment hull is the closed convex hull of the pairs
$(z,zz^{\top})$ with $z\in\Z^n$. We prove that one round of
Chv\'atal--Gomory (CG) cuts applied to its basic semidefinite relaxation
recovers this hull exactly if and only if $n\le6$. Our approach identifies
its facets with perfect Delaunay polyhedra and characterizes their recovery
through depth, the gap between lattice and continuous minima of a defining
quadratic after normalization. A facet inequality is valid for the
rank-one CG closure exactly when its depth is less than one. We also
analyze the complexity of deciding
membership in the CG closure and finding violated inequalities. Every
rational point outside the CG closure admits a violated CG inequality whose
encoding length is polynomial in the input length. Membership is strongly
coNP-complete and exact separation is strongly NP-hard. The same reduction
proves strong NP-hardness of quadratic split-inequality separation,
resolving a conjecture of Burer and Letchford.
\end{abstract}



\input{sections/01_introduction}
\input{sections/02_framework}
\input{sections/03_facets}
\input{sections/04_depth}
\input{sections/05_exact_hull}
\input{sections/06_complexity}
\input{sections/07_conclusion}
\input{appendices/A_low_dimensional}
\input{appendices/B_encoding_bounds}

\section*{Declaration of AI use}
The authors used GPT to assist with reviewing the manuscript, refining its
exposition, and preparing the bibliography. They reviewed all AI-assisted
material and take full responsibility for the manuscript.

\bibliographystyle{unsrtnat-names}
\bibliography{references}
\end{document}

%% file: preamble.tex
\usepackage[T1]{fontenc}
\usepackage{lmodern}
\usepackage{microtype}
\usepackage{hyphenat}
\usepackage{amsmath,amssymb,amsthm,mathtools}
\usepackage{booktabs}
\usepackage[normalem]{ulem}
\usepackage{xcolor}
\usepackage{enumitem}
\usepackage[margin=1in]{geometry}
\usepackage[numbers,sort&compress]{natbib}
\usepackage[colorlinks=true,linkcolor=blue!55!black,
            citecolor=blue!55!black,urlcolor=blue!55!black]{hyperref}

\numberwithin{equation}{section}

\newtheorem{theorem}{Theorem}[section]
\newtheorem{proposition}[theorem]{Proposition}
\newtheorem{definition}[theorem]{Definition}
\newtheorem{lemma}[theorem]{Lemma}
\newtheorem{corollary}[theorem]{Corollary}
\theoremstyle{plain}

\newcommand{\R}{\mathbb R}
\newcommand{\Z}{\mathbb Z}
\newcommand{\Q}{\mathbb Q}
\newcommand{\Smat}{\mathbb S}
\newcommand{\MIQ}{\operatorname{MIQ}}
\newcommand{\Erd}{\operatorname{Erdahl}}
\newcommand{\conv}{\operatorname{conv}}
\newcommand{\aff}{\operatorname{aff}}
\newcommand{\lin}{\operatorname{span}}

\newcommand{\range}{\operatorname{range}}
\newcommand{\rank}{\operatorname{rank}}
\newcommand{\relint}{\operatorname{relint}}
\newcommand{\GL}{\operatorname{GL}}
\newcommand{\CG}{\operatorname{CG}}
\newcommand{\Space}{\operatorname{Space}}

\newcommand{\ip}[2]{\left\langle #1,#2\right\rangle}

\newcommand{\CGViol}{\textnormal{\textsc{CG-Violation}}}
\newcommand{\CGMem}{\textnormal{\textsc{CG-Membership}}}
\newcommand{\F}{\mathbb F}

%% file: sections/01_introduction.tex
\section{Introduction}
\label{sec:introduction}

Integer quadratic optimization is an important class of nonlinear integer programming problems. Consider any quadratic function $q:\R^n \rightarrow\R$, 
\[q(z)=z^\top A z+2b^\top z+\alpha,\]
where $A\in \R^{n\times n}$ and is symmetric, $b\in\mathbb R^n$, and $\alpha\in\mathbb R$. In this work, we focus on optimization of $q$ over $z\in \Z^n$, namely $\inf_{z\in\Z^n} q(z)$ or $\sup_{z\in\Z^n} q(z)$. We impose no convexity or concavity assumption on $q$. We use moment lifting to represent each quadratic objective as a linear functional of a matrix containing the first- and second-order terms. The closed convex hull of the lifted integer points then provides a common convex set for both minimization and maximization.

\subsection{The moment hull and its relaxations}

Throughout, we use $\R$, $\Q$, and $\Z$ to denote the real numbers, rational numbers, and integers, respectively. Let $\Smat^n$ denote the space of real symmetric $n\times n$ matrices for every positive $n\in \Z$, and let $\langle A,B\rangle=\operatorname{tr}(AB)$ denote the trace inner product. For any $z\in\mathbb R^n$, its lifted moment matrix is defined by
\[
W_z :=
\begin{bmatrix}
1\\
z
\end{bmatrix}
\begin{bmatrix}
1\\
z
\end{bmatrix}^{\!\top}
=
\begin{bmatrix}1&z^{\top}\\ z&zz^{\top}\end{bmatrix}.
\]
For the arbitrary quadratic function $q(z)=z^\top A z+2b^\top z+\alpha$, let
\[
M_q :=
\begin{bmatrix}
\alpha & b^\top\\
b & A
\end{bmatrix},
\]
then
\[
q(z)=\langle M_q,W_z\rangle. 
\]
Thus, moment lifting equivalently expresses a quadratic function over $z\in \R^n$ by a
linear function over lifted matrix $W_z$.

Specifically, when $z$ is restricted to $\Z^n$, moment lifting naturally gives rise to the set of all lifted matrices associated with $z\in \Z^n$. We refer to the closure of its convex hull by \emph{integer quadratic moment hull} and denote it by 
\[
\mathrm{MIQ}_n
:=
\operatorname{cl}\operatorname{conv}
\{W_z:z\in\mathbb Z^n\}
\]
throughout this work. Our notation $\MIQ_n$ corresponds to the purely integer case $\mathrm{MIQ}_{n,0}$ of \citet{BurerLetchford2014}.

Throughout, we index the rows and columns of an $(n+1)\times(n+1)$ moment matrix by $\{0,1,\ldots,n\}$. We let 
\[
\mathcal{H}_n:=\{W\in\Smat^{n+1}: W_{00}=1 \}, 
\]
and denote the semidefinite relaxation of $\MIQ_n$ by
\[
\mathcal P_n:=\{W\in\mathcal H_n: W\succeq 0\}.
\]
Every lattice moment matrix $W_z$ is positive semidefinite, so
\[
\mathrm{MIQ}_n\subseteq\mathcal P_n\subseteq\mathcal H_n.
\]

Our focus is the \emph{rank-one Chv\'atal--Gomory (CG) closure}
$\mathcal C_n$ of $\mathcal P_n$, obtained by applying the CG procedure once
in the natural integral moment coordinates. To assess the strength of CG
rounding, we compare $\mathcal C_n$ with the \emph{split-inequality closure}
$\mathcal S_n$, defined by the quadratic split inequalities introduced by
Burer and Letchford~\citep{BurerLetchford2014}.
Section~\ref{sec:rank-one-closure} gives the formal definitions. These sets satisfy
\[
 \MIQ_n\subseteq\mathcal C_n\subseteq\mathcal S_n\subsetneq\mathcal P_n.
\]
This chain gives a baseline comparison. However, we would like to further explore the following research questions. First, how far is $\mathcal C_n$ from fully recovering $\MIQ_n$, and how does $\mathcal C_n$ compare with $\mathcal S_n$? Second, how hard is exact separation over
$\mathcal C_n$ when the dimension is part of the input? We first review the relevant literature.

\subsection{Related literature}

Quadratic optimization with integer decisions arises in wide-ranging applications, such as portfolio selection and facility location \cite{GunlukLinderoth2010}, unit commitment \cite{FrangioniGentile2006}, sparse regression \cite{WeiGomezKucukyavuz2022}, and learning Gaussian Bayesian networks \cite{KucukyavuzEtAl2023Bayesian}. The breadth of these applications has sustained interest in the geometry of quadratic integer optimization, with convex-hull descriptions and valid inequalities providing a natural bridge from discrete structure to stronger continuous relaxations.

Several earlier works study quadratic optimization by lifting first and second moments. For binary variables, Padberg's Boolean quadric polytope provides the basic polyhedral model for standard quadratic linearizations
\cite{Padberg1989}.  For unconstrained quadratic $0$--$1$ minimization, Boros, Crama, and Hammer
gave a polynomial-size linear programming formulation of the first
Chv\'atal closure of the standard linearization relaxation. They also
related this closure to odd-cycle inequalities through a signed-graph
reformulation \cite{BorosCramaHammer1992}.  Burer gave an exact formulation over the completely positive cone for a broad class of nonconvex quadratic programs with binary and nonnegative continuous variables
\cite{Burer2009}.  Anstreicher and Burer obtained exact
representations of moment convex hulls over several classes of low-dimensional
polytopes \cite{AnstreicherBurer2010}. Moreover,
Burer and Letchford developed a moment-hull framework for mixed-integer
quadratic programming \cite{BurerLetchford2014}.

For the purely integer moment hull $\MIQ_n$, Burer and Letchford proved that
their quadratic split
inequalities dominate the positive-semidefinite inequalities
\cite[Theorem~4]{BurerLetchford2014} and completely describe the moment hull in
dimension two: $\MIQ_2=\mathcal S_2$
\cite[Theorem~8]{BurerLetchford2014}. In particular, they left open whether split inequalities suffice in any dimension above two. In an attempt to address this, in dimension six, they provided a bounded facet whose defining inequality is not a split inequality and concluded that split inequalities do not completely describe the hull \cite[Section~5.3]{BurerLetchford2014}. The authors also suggested a connection between facets of the integer quadratic moment hull and the structures called Delaunay polytopes. We include relevant background information on them in Section~\ref{sec:delaunay}. Our geometric approach in this work is inspired by this connection. Buchheim and Traversi later recorded that exactness still remained unresolved in dimensions three through five \cite[Section~3, discussion following Theorem~4]{BuchheimTraversi2015}. They gave an exact separation method for these inequalities and evaluated it computationally. For positive-semidefinite input, the separation reduces to convex quadratic integer minimization. In addition, for input outside the positive-semidefinite cone, their Algorithm~1 finds a violated split inequality in polynomial time \cite[Section~4, Theorem~5]{BuchheimTraversi2015}.

The quadratic split inequalities discussed above use the same lattice
disjunctions as classical split cuts for mixed-integer linear programs. Cook,
Kannan, and Schrijver introduced split cuts for rational mixed-integer
polyhedra and proved that the resulting split closure is polyhedral
\cite{CookKannanSchrijver1990}.  Andersen,
Cornu\'ejols, and Li showed that intersection cuts completely describe the
split closure, provided infeasible and nonoptimal bases are included
\cite{AndersenCornuejolsLi2005}.  \c{C}ezik and
Iyengar extended Chv\'atal--Gomory, lift-and-project, and related cutting-plane
constructions to mixed $0$--$1$ conic programs \cite{CezikIyengar2005}. 

Our analysis of $\mathcal C_n$ also draws on the broader theory of
CG closures. CG procedure builds on the foundational works by
Gomory and Chv\'atal \cite{Gomory1958,Chvatal1973}. Based on the result that the elementary closure of a rational polyhedron is again a rational polyhedron by \citet{Schrijver1980}, Dadush, Dey, and Vielma extended this to strictly convex bodies and their intersections with rational polyhedra \cite{DadushDeyVielma2011}.  They subsequently proved that the CG closure of every nonempty compact convex set is a rational polytope \cite{DadushDeyVielma2014}.  On the other hand, they exhibited a strictly convex body in $\R^3$ whose split closure is not polyhedral \cite{DadushDeyVielmaSplit2011}.
De Meijer and Sotirov systematically developed the CG procedure for integer
semidefinite programs and spectrahedra \cite{deMeijerSotirov2024}.  In contrast, we study a highly structured but unbounded moment spectrahedron.

By exploiting specific constraint and objective structures, constrained mixed- and pure-integer variants of quadratic optimization have also been studied through strong formulations and valid inequalities. For models in which binary indicators switch continuous variables on or off, perspective cuts and reformulations strengthen convex quadratic relaxations
\cite{FrangioniGentile2006,GunlukLinderoth2010}. Perspective strengthening has also demonstrated success combined with outer approximation in mixed-integer convex formulations for learning Gaussian directed acyclic graphs
\cite{XuEtAl2025DAG}. Moreover, semidefinite methods yield equivalent convex quadratic reformulations for classes of constrained models
\cite{billionnet2012extending} and exact extended formulations for convex quadratic epigraphs with indicators \cite{WeiAtamturkGomezKucukyavuz2024}. Special structures in the quadratic coefficient matrix give rise to further valid inequalities \cite{AtamturkGomez2018,AtamturkGomez2023}. A different approach uses coordinate-wise optimality conditions to derive optimality-preserving reformulations for convex quadratic programs with indicators \cite{xu2026coordinate}. Extensions of submodularity are applied more broadly
to structured nonlinear minimization over constrained mixed- and
pure-integer feasible sets, yielding epigraph convex-hull descriptions and
exact separation algorithms \cite{YuKucukyavuz2025DR,yu2025convexification}. These results cover special cases of quadratic integer optimization. In our work, we focus on quadratic optimization over the full integer lattice, without imposing structural restrictions on the quadratic objective.

A substantial literature addresses the complexity of elementary closures and
cut separation.  Eisenbrand proved that deciding whether a given rational
point lies outside the elementary closure is NP-complete
\cite{Eisenbrand1999}.  Bockmayr and Eisenbrand proved that, in fixed
dimension, the elementary closure of a rational polyhedron can be described by
polynomially many inequalities and constructed in polynomial time
\cite{BockmayrEisenbrand2001}.  Cornu\'ejols and Li compared the elementary
closures generated by several classical families of cutting planes
\cite{CornuejolsLi2001}.

Caprara and Fischetti proved that the separation problem for
$\{0,\tfrac12\}$-CG cuts is NP-complete \cite{CapraraFischetti1996}.  Caprara
and Letchford established strong NP-completeness for split-cut separation and
related families \cite{CapraraLetchford2003}.
Fischetti and Lodi modeled rank-one CG separation as a mixed-integer program
and used the resulting separator in a computational study of optimization over
the first Chv\'atal closure \cite{FischettiLodi2007}.  Dash,
G\"unl\"uk, and Lodi proved that the MIR closure coincides with the split
closure and that its separation problem is NP-hard \cite{DashGunlukLodi2010}.
For the quadratic split inequalities considered here, Burer and Letchford
conjectured that separation is strongly NP-hard
\cite[p.~255]{BurerLetchford2014}.

\subsection{Our contributions}

Our results identify the exact dimension threshold for recovery of the integer
quadratic moment hull by one round of Chv\'atal--Gomory cuts. We also settle
the open question of whether quadratic split inequalities describe the moment
hull in dimensions three through five
\cite{BurerLetchford2014,BuchheimTraversi2015} and prove Burer and Letchford's
conjecture that separating these inequalities is strongly NP-hard
\cite[p.~255]{BurerLetchford2014}.

For every $n\ge1$, the rank-one CG closure of the basic semidefinite
relaxation $\mathcal P_n$ satisfies
\[
 \mathcal C_n=\MIQ_n
 \quad\Longleftrightarrow\quad
 n\le6.
\]
The split-inequality closure is exact only through dimension five, so one CG
round extends exact recovery by one dimension. The complete comparison is
\[
\begin{array}{rll}
1\le n\le5: & \MIQ_n=\mathcal C_n=\mathcal S_n, &\\[1mm]
n=6:         & \MIQ_n=\mathcal C_n\subsetneq\mathcal S_n, &\\[1mm]
n\ge7:       & \MIQ_n\subsetneq\mathcal C_n\subsetneq\mathcal S_n. &
\end{array}
\]
Moreover, $\mathcal S_n\subsetneq\mathcal P_n$ in every dimension.

Our geometric results explain these thresholds. We prove a bijection between
the facets of $\MIQ_n$ and perfect Delaunay polyhedra, making precise the
connection suggested by Burer and Letchford
\cite[Section~5.3]{BurerLetchford2014}. We associate a depth with each facet:
the gap between the lattice and continuous minima of a primitive
integer-valued defining quadratic. We prove that the facet is valid for
$\mathcal C_n$ if and only if its depth is less than one. Exact depth
calculations for the known types through dimension seven, together with
extension of their defining quadratics to higher dimensions, give the
comparisons above.

We also establish complexity results for membership and exact separation over
$\mathcal C_n$. Every rational point outside $\mathcal C_n$ admits a violated
rank-one CG inequality whose encoding length is polynomial in the input
length. When $n$ is part of the input, nonmembership in $\mathcal C_n$ is
strongly NP-complete and membership is strongly coNP-complete, even when the
input data are polynomially bounded in $n$. The same reduction proves strong
NP-hardness of quadratic split-inequality separation, even for input points
in $\mathcal P_n$.

\subsection{Organization}

Section~\ref{sec:framework} connects quadratic functions with linear
inequalities in the moment variables. It introduces the Erdahl cone,
consisting of quadratics that are nonnegative at every integer point, and
Delaunay polyhedra, the nonempty sets of integer points where a nonzero such
quadratic vanishes. It also specifies the coordinates used for CG rounding and defines
the CG and split-inequality closures.
Section~\ref{sec:facets} proves the correspondence between the facets of
$\MIQ_n$ and perfect Delaunay polyhedra.
Section~\ref{sec:depth} determines which facet inequalities hold after one
round of CG cuts. The key quantity is depth, which measures the gap between
the integer and continuous minima of a defining quadratic after choosing a
standard scale for its coefficients. We express CG inequalities using
integer-valued quadratics, called cut polynomials, and prove that a facet
inequality is valid for $\mathcal C_n$ exactly when its depth is less than
one. We further provide a  rigidity argument to establish that a facet inequality failing this
test cannot follow from the other cuts in the same round.
Section~\ref{sec:low-dimensional} applies known classifications of perfect
Delaunay polyhedra through dimension seven and computes their depths. By
extending suitable examples to higher dimensions, we determine exactly when
the CG- and split-inequality closures recover the integer quadratic moment
hull. In Section~\ref{sec:complexity}, we study membership and separation for
$\mathcal C_n$. Specifically, we show that every rational point outside the closure has a
short certificate of nonmembership that can be checked efficiently. We then
establish the computational difficulty of deciding membership and finding
violated inequalities when the dimension is part of the input.
Section~\ref{sec:conclusion} summarizes the results and discusses directions
for further research.

%% file: sections/02_framework.tex
\section{Quadratic and geometric foundations}
\label{sec:framework}

This section provides the algebraic and geometric framework used throughout this paper. In Section \ref{sec:quadratics_erdahl}, we identify the correspondence between the valid affine inequalities for $\MIQ_n$ and the elements of the so-called \emph{Erdahl cone}. Section \ref{sec:delaunay} introduces Delaunay polyhedra and the notion of perfection, which are later used to describe the facets of $\MIQ_n$. Then, in Section~\ref{sec:coordinates}, we discuss the integral moment coordinates and relate them to the integer-valued quadratic polynomials. Lastly, Section \ref{sec:rank-one-closure} formally defines the CG closure, $\mathcal C_n$, and the split-inequality closure, $\mathcal S_n$. We further prove their containment relation and remark on their invariance under integer-affine changes of variables. 
We use $\conv$, $\aff$, and $\range$ to denote convex hull, affine hull,
and range, respectively.

\subsection{Quadratic polynomials and the Erdahl cone}
\label{sec:quadratics_erdahl}

Let $\mathcal Q_n$ be the vector space of real polynomials in $n$ variables of degree at most two. Every $q\in\mathcal Q_n$ has a unique representation
\[
 q(z) = z^{\top}A_qz + 2b_q^{\top}z+\alpha_q,
 \]
 where $A_q\in\Smat^n$, $b_q\in\R^n$, and $\alpha_q\in\R$. Equivalently, with
 \[
 M_q:=\begin{bmatrix}\alpha_q&b_q^{\top}\\ b_q&A_q\end{bmatrix}
 \in\Smat^{n+1},
\]
for $W=\left[\begin{smallmatrix}1&x^{\top}\\x&X\end{smallmatrix}\right]
\in\mathcal H_n$, define the linearization
\begin{equation}
\label{eq:linearization}
 L_q(W) := \ip{M_q}{W} =\ip{A_q}{X}+2b_q^{\top}x+\alpha_q.
\end{equation}
Then $L_q(W_z) = q(z)$ for every $z\in\R^n$. 
The relevant dimensions follow from the number of independent matrix entries. A symmetric $(n+1)\times(n+1)$ matrix is determined by its $n+1$ diagonal entries and its $\binom{n+1}{2}$ entries above the diagonal. Hence,
\[
\dim \mathbb S^{n+1}
=
(n+1)+\binom{n+1}{2}
=
\binom{n+2}{2}.
\]
The constraint $W_{00}=1$ removes one degree of freedom, so
\[
\dim \mathcal H_n
=
\binom{n+2}{2}-1.
\]
For brevity, we will consistently use 
\[
d_n
:=
\binom{n+2}{2}-1
=
n+\binom{n+1}{2}
\]
to represent this constant. Similarly, $\dim \mathcal Q_n=d_n+1$.

The \emph{Erdahl cone} is defined by
\[
 \Erd(n):=\{q\in\mathcal Q_n:q(z)\ge 0,\, \forall\, z\in\Z^n\}.
\]
We call every $q\in\Erd(n)$ a \emph{lattice-nonnegative quadratic} to highlight the property that $q$ is nonnegative at all integer arguments. We remark that such a quadratic need not be nonnegative on all of $\R^n$, nor need it take integer values on $\Z^n$.

We first record an elementary property of quadratic polynomials that will
be used to identify affine functionals from their values at the lattice
moment matrices.

\begin{lemma}[Quadratics vanishing on the integer lattice]
\label{lem:lattice-zero-polynomial}
If $p\in\mathcal Q_n$ satisfies $p(z)=0$ for every $z\in\Z^n$, then
$p\equiv0$.
\end{lemma}

\begin{proof}
By writing $p(z)=z^{\top}A_pz+2b_p^{\top}z+\alpha_p$ and evaluating it at $z=0$, we obtain $\alpha_p = 0$. Let $e_i$ be the $i$th standard unit vector. Evaluating
at $e_i$ and $-e_i$ gives
\[
 (A_p)_{ii}+2(b_p)_i=0,
 \qquad
 (A_p)_{ii}-2(b_p)_i=0,
\]
so $(A_p)_{ii}=(b_p)_i=0$ for every $i$. For $i<j$, evaluating at
$e_i+e_j$ then gives $2(A_p)_{ij}=0$. Thus $A_p=0$, $b_p=0$, and
$\alpha_p=0$, proving that $p\equiv0$.
\end{proof}

An immediate consequence is that the lattice moment matrices $W_z$,
$z\in\Z^n$, affinely span $\mathcal H_n$. Indeed, every affine functional
$F$ on $\mathcal H_n$ can be written as $F=L_p$ for some
$p\in\mathcal Q_n$. If $F(W_z)=0$ for every $z\in\Z^n$, then
$p(z)=0$ for every $z\in\Z^n$, so
Lemma~\ref{lem:lattice-zero-polynomial} gives $p\equiv0$ and hence
$F\equiv0$. Consequently, $\MIQ_n$ is full-dimensional in
$\mathcal H_n$ and has dimension $d_n$.

For $q\in\mathcal Q_n$, the validity characterization of
\citet{BurerLetchford2014} is
\begin{equation}\label{eq:moment-validity}
 L_q(W)\ge0\text{ for every }W\in\MIQ_n
 \quad\Longleftrightarrow\quad
 q\in\Erd(n).
\end{equation}
Validity implies $q(z)=L_q(W_z)\ge0$ at every integer point.
Conversely, if $q\in\Erd(n)$, then $L_q$ is nonnegative at every lattice
moment matrix and hence on their closed convex hull.

\subsection{Delaunay polyhedra and perfection}
\label{sec:delaunay}

For $q\in\Erd(n)$, define its \emph{lattice-zero set} by
\[
 Z(q):=\{z\in\Z^n:q(z)=0\}.
\]
\begin{definition}[Delaunay polyhedron]
\label{def:delaunay_polyhedron}
A nonempty set $D\subseteq\Z^n$ is a \emph{Delaunay polyhedron} if
$D=Z(q)$ for some nonzero $q\in\Erd(n)$.
\end{definition}

Throughout the paper, $D$ denotes the lattice-point set and $\conv(D)$ denotes
the associated convex set.  We call $D$ bounded when $\conv(D)$ is bounded,
or equivalently when $D$ is finite.  In this case, $\conv(D)$ is the associated
Delaunay polytope.  When no confusion can arise, we also call the finite
lattice-point set $D$ a Delaunay polytope.

To formalize when a Delaunay polyhedron uniquely determines its defining
quadratic up to scaling, for an arbitrary $D\subseteq\Z^n$, let
\[
 \Space(D):=\{p\in\mathcal Q_n:p(z)=0\text{ for every }z\in D\}.
\]

\begin{definition}[Perfect Delaunay polyhedron] 
\label{def:perfect_delaunay}
A Delaunay polyhedron $D$ is \emph{perfect} if
\[
 \dim\Space(D)=1.
\]
\end{definition}

If $D$ is perfect, $\Space(D)\cap\Erd(n)$ is a nonzero ray. We denote
a generator of this ray by $q_D$, which is unique up to multiplication by a
positive scalar.

These notions arise in the theory of inhomogeneous quadratic forms. Erdahl studied the cone of quadratics that are nonnegative on $\Z^n$, characterized
its extreme rays, and connected one of the extreme-ray classes to perfect
Delaunay polyhedra through lattice-zero sets \cite{Erdahl1992}. For further
background on lattice geometry and computational methods, see Sch\"urmann
\cite{Schuermann2008} and Dutour Sikiri\'c, Sch\"urmann, and Vallentin
\cite{DutourSikiricSchuermannVallentin2009}.

In fact, every perfect Delaunay polyhedron is full-dimensional in $\R^n$.  Suppose otherwise that $\aff(D)\subsetneq \R^n$. Then there exists a nonzero affine polynomial $q\in\mathcal Q_n$ that vanishes on $\aff(D)$; that is, $q(z) = 0$ for every $z\in \aff(D)$. 
To see this, we write $\aff(D)$ as $z^0 + V$, where $z^0\in D$ and $V$ is a linear subspace with $\dim V\leq n-1$. Its orthogonal complement
\[V^\perp := \{u\in\R^n : u^\top v = 0,\,\forall\, v\in V\}\]
has dimension $n - \dim V \geq 1$. Choose $u\in V^\perp$ with $u\neq 0$, then $q(z) := u^\top(z - z^0)$ is a nonzero affine polynomial in $\mathcal Q_n$ that vanishes on $\aff(D)$.
Hence, $q\in \Space(D)$. Meanwhile, $q_D$ cannot be affine because a nonconstant affine polynomial is unbounded below on $\Z^n$, and a nonzero constant has no lattice zeros. Therefore, $q$ and $q_D$ are linearly independent, contradicting $\dim\Space(D)=1$.

Perfect Delaunay polyhedra may nevertheless be unbounded.
Section~\ref{subsec:bounded-facets} describes their structure and its relation
to facet boundedness.

\subsection{Matrix coordinates and integral moment coordinates}
\label{sec:coordinates}
For $W=\left[\begin{smallmatrix}1&x^{\top}\\x&X\end{smallmatrix}\right]
\in\mathcal H_n$, we call $(x,X)$ its \emph{matrix coordinates}, with $X$
represented by its independent upper-triangular entries.

The CG procedure uses affine functionals that are integer-valued at every lattice moment point.  In monomial coordinates, however, an integer-valued quadratic may have fractional coefficients. For example, let \[\binom{z}{2} := z(z-1)/2\]
for $z\in\R$; here, $\binom{z}{2}$ is integer-valued on $z\in\Z$.  We therefore use the Newton basis of $\mathcal Q_n$:
\[
 1,\qquad z_i,\qquad \binom{z_i}{2},\qquad z_iz_j,
\]
for $j\in \{i+1, \dots, n\}$, $i\in \{1,\dots,n\}$.
This basis gives the natural integral coordinates for CG rounding.

The \emph{integral moment coordinates} of $W$ are
\[
 \eta(W)=
 \left((x_i)_{i=1}^n;\;
       (w_{ii})_{i=1}^n;\;
       (w_{ij})_{1\le i<j\le n}\right) = (x,w) \in\R^{d_n}. 
\]
Here, $w$ is given by
\begin{equation}\label{eq:eta-coordinates}
 w_{ij} := \begin{cases}
 (X_{ii}-x_i)/2, & \forall\, 1\leq i = j \leq n, \\
 X_{ij}, & \forall\, 1\le i<j\le n. 
 \end{cases}
\end{equation}
The diagonal coordinates $w_{11},\ldots,w_{nn}$ precede the off-diagonal
coordinates $w_{ij}$, whose index pairs $(i,j)$ with $i<j$ are ordered
lexicographically.

The coordinate lattice for CG rounding is $\Z^{d_n}$ in the $(x,w)$
representation. In matrix coordinates, the corresponding affine lattice
consists of matrices satisfying $x_i\in\Z$ and $X_{ii}-x_i\in2\Z$ for
every $i$, and $X_{ij}\in\Z$ for every $i<j$. We denote the inverse coordinate map by $W:\R^{d_n}\rightarrow\mathcal H_n$.
For $(x,w)\in\R^{d_n}$, define
\[
 W(x,w):=
 \begin{bmatrix}1&x^{\top}\\x&X(x,w)\end{bmatrix},
\]
where $X(x,w)\in\Smat^n$ is given by
\[X(x,w)_{ij} := 
\begin{cases}
    x_i + 2w_{ii}, & \forall\, 1\leq i = j\leq n, \\
    w_{ij}, & \forall\, 1\leq i < j \leq n, \\
    w_{ji}, & \forall\, 1\leq j < i \leq n. 
\end{cases}
\]

For any rank-one moment matrix $W_z$ associated with $z\in \R^n$, we define 
\begin{equation}\label{eq:eta-generator}
\eta(z) := \eta(W_z)=
 \left((z_i)_{i=1}^n;\;
       \left(\binom{z_i}{2}\right)_{i=1}^n;\;
       (z_iz_j)_{1\le i<j\le n}\right)
\end{equation}
by slight abuse of notation. In particular, $\eta(z)\in\Z^{d_n}$ for every $z\in\Z^n$.

Every $p\in\mathcal Q_n$ has a unique representation in this Newton basis: 
\[
 p(z)=\alpha_p+\sum_{i=1}^n a_i z_i
 +\sum_{i=1}^n d_i\binom{z_i}{2}
 +\sum_{1\le i<j\le n}\gamma_{ij}z_iz_j
 =\alpha_p+c_p^{\top}\eta(z),
\]
where
\[
 c_p=\left((a_i)_{i=1}^n;\;
           (d_i)_{i=1}^n;\;
           (\gamma_{ij})_{1\le i<j\le n}\right)
 \in\R^{d_n}.
\]
We call $c_p$ the \emph{integral moment normal} of $p$, and write
\[
 \operatorname{Int}_2(\Z^n)
 :=\{p\in\mathcal Q_n:p(\Z^n)\subseteq\Z\}.
\]

The following lemma verifies the correspondence between the two coordinate
systems and shows that a quadratic is integer-valued on $\Z^n$ exactly when
its Newton coefficients are integers.

\begin{lemma}[Integral moment coordinates and integer-valued quadratics] 
\label{lem:integral-moment-coordinates}
\mbox{}
\begin{itemize}
    \item[(i)] The mapping $\eta:\mathcal H_n \rightarrow\R^{d_n}$ given in \eqref{eq:eta-coordinates} is an affine bijection, and $\eta(W_z)=\eta(z)\in\Z^{d_n}$ for every $z\in\Z^n$.
    \item[(ii)] For every $p\in\mathcal Q_n$ and every $W\in \mathcal H_n$, 
\begin{equation}\label{eq:evaluation-dictionary}
 L_p(W)=\alpha_p+c_p^{\top}\eta(W).
\end{equation}
    \item[(iii)] A polynomial $p\in \mathcal Q_n$ satisfies  $p\in\operatorname{Int}_2(\Z^n)$ if and only if $(\alpha_p,c_p)\in\Z\times\Z^{d_n}$. 
\end{itemize}
\end{lemma}

\begin{proof}
For (i), the formulas above show directly that $\eta$ and $W$ are inverse affine maps by construction. Substituting $x=z$ and $X=zz^{\top}$ gives \eqref{eq:eta-generator}, whose entries are integral when $z\in \Z^n$.

For (ii), the two sides of \eqref{eq:evaluation-dictionary} agree at every
$W_z$. They are affine functions on $\mathcal H_n$, so the affine-spanning
consequence of Lemma~\ref{lem:lattice-zero-polynomial} implies that they
agree everywhere.

Lastly, we prove claim (iii). If $p$ is integer-valued on $\Z^n$, its Newton coefficients are recovered by
\[
 \alpha_p=p(0),\qquad
 a_i=p(e_i)-p(0),\qquad
 d_i=p(2e_i)-2p(e_i)+p(0),
\]
for $1\leq i\leq n$, and for $1\leq i<j \leq n$,
\[
 \gamma_{ij}=p(e_i+e_j)-p(e_i)-p(e_j)+p(0).
\]
Thus all Newton coefficients are integers. Conversely, every Newton-basis element is integer-valued on $\Z^n$, so every integral linear combination of these functions is integer-valued. 
\end{proof}

The Newton representation is convenient for characterizing integer-valued quadratics, while the standard matrix and Hessian representations will be useful for analyzing convexity and continuous minima. We therefore summarize the correspondence between these representations:

\begin{equation}\label{eq:coefficient-dictionary}
 p(z)=z^{\top}A_pz+2b_p^{\top}z+\alpha_p
     =\frac12z^{\top}H_pz+\ell_p^{\top}z+\alpha_p,
\end{equation}
where
\begin{equation}\label{eq:newton-to-matrix}
 (A_p)_{ii}=\frac{d_i}{2},\qquad
 (A_p)_{ij}=\frac{\gamma_{ij}}{2}\ (i<j),\qquad
 (b_p)_i=\frac{a_i}{2}-\frac{d_i}{4},
\end{equation}
and
\begin{equation}\label{eq:newton-to-hessian}
 H_p=2A_p,\qquad
 (H_p)_{ii}=d_i,\qquad
 (H_p)_{ij}=\gamma_{ij}\ (i<j),\qquad
 (\ell_p)_i=a_i-\frac{d_i}{2}.
\end{equation}
Thus $H_p$ is integral whenever $p$ is integer-valued, although $A_p$, $b_p$,
and $\ell_p$ need not be integral.

In what follows, under the affine bijection $\eta$, we will use the matrices and coordinate vectors interchangeably, while keeping the same notation for sets and functions in either representation.  For example, in integral moment coordinates, the semidefinite relaxation is written as
\[
 \mathcal P_n=
 \{\eta\in\R^{d_n}:W(\eta)\succeq0\}.
\]

\subsection{The rank-one CG closure and quadratic split inequalities}
\label{sec:rank-one-closure}
This section defines the rank-one CG closure $\mathcal C_n$ and the
split-inequality closure $\mathcal S_n$. We establish their basic
containment relations with $\mathrm{MIQ}_n$ and $\mathcal P_n$, and show
that all these four sets are invariant under integer-affine bijections.

For a nonempty closed convex set $K\subseteq\R^{d_n}$ and
$c\in\Z^{d_n}$, let
\[
 \beta_K(c):=\inf\{c^{\top}\eta:\eta\in K\}.
\]
When $\beta_K(c)>-\infty$, the associated CG cut is $c^\top\eta \geq \lceil\beta_K(c)\rceil$. The CG closure of $K$ with respect to the integral moment lattice is
\[
 \CG(K):=
 \bigcap_{\substack{c\in\Z^{d_n}\\\beta_K(c)>-\infty}}
 \left\{\eta:c^{\top}\eta\ge\lceil\beta_K(c)\rceil\right\}.
\]
We now specialize this construction to the semidefinite relaxation
$\mathcal P_n$. For $c\in\Z^{d_n}$, write the integer-valued quadratic
$c^{\top}\eta(x)$, $x\in\R^n$, as
\[
 c^{\top}\eta(x)=x^\top A_cx+2b_c^\top x.
\]
Then
\begin{equation}\label{eq:beta-formula}
 \beta_{\mathcal P_n}(c)=\inf_{x\in\R^n}c^{\top}\eta(x)
 =\begin{cases}
 -b_c^\top A_c^\dagger b_c,
& A_c\succeq 0\text{ and }b_c\in\range(A_c),\\
 -\infty,&\text{otherwise}.
 \end{cases}
\end{equation}
Here $A_c^\dagger$ denotes the Moore-Penrose pseudoinverse.  To verify the
formula, write $X=xx^{\top}+S$ with $S\succeq0$.  The linearized objective is
\[
 c^{\top}\eta(x)+\ip{A_c}{S}.
\]
If $A_c\not\succeq0$, the objective is unbounded below as $S$ varies,
and $c^{\top}\eta(x)$ is unbounded below along a direction of negative
curvature. If $A_c\succeq0$, the second term is nonnegative and is minimized
at $S=0$.
The problem then reduces to minimizing $c^{\top}\eta(x)$ over $\R^n$.  Its value
is finite exactly when $b_c\in\range(A_c)$, which gives
\eqref{eq:beta-formula}.  In that case, $x=-A_c^{\dagger}b_c$ is a minimizer,
and the corresponding moment matrix $W_x$ has matrix rank one.

With this characterization of the relaxation value, we define the \emph{rank-one CG closure} of $\mathcal P_n$  by
\begin{equation}\label{eq:Cn}
 \mathcal C_n:=\CG(\mathcal P_n)
 =\bigcap_{\substack{c\in\Z^{d_n}\\\beta_{\mathcal P_n}(c)>-\infty}}
 \left\{\eta:c^{\top}\eta
 \ge\left\lceil\beta_{\mathcal P_n}(c)\right\rceil\right\}.
\end{equation}
Here, ``rank-one'' means that the CG procedure is applied exactly once to
$\mathcal P_n$.
This is the elementary CG closure of the spectrahedron $\mathcal P_n$ in the
integral coordinates \eqref{eq:eta-coordinates} (see
\cite[Definition~1 and Sections~2.1--2.2]{deMeijerSotirov2024}).

Although \eqref{eq:Cn} is written solely as an intersection of affine
halfspaces, the semidefinite constraint defining $\mathcal P_n$ is
already implied by this family of CG inequalities. We now verify this fact. 

For
$v=(v_0,v')\in\Z^{n+1}$, consider the integer-valued affine square
\[
 q_v(z):=\left(v_0+{v'}^{\top}z\right)^2.
\]
Write $q_v(z)=v_0^2+c_v^{\top}\eta(z)$, where $c_v\in\Z^{d_n}$.
Since $q_v\ge0$ on $\R^n$, \eqref{eq:beta-formula} gives
$\beta_{\mathcal P_n}(c_v)\ge-v_0^2$.  The CG inequality for $c_v$ therefore
gives, for every $\eta\in\mathcal C_n$,
\[
 c_v^{\top}\eta\ge\left\lceil\beta_{\mathcal P_n}(c_v)\right\rceil
 \ge-v_0^2,
\]
because $v_0^2$ is an integer.  Hence
$v^{\top}W(\eta)v=v_0^2+c_v^{\top}\eta\ge0$ for every
$\eta\in\mathcal C_n$ and every $v\in\Z^{n+1}$.

These integral inequalities characterize positive semidefiniteness.  Indeed,
if $W\not\succeq0$, the open set
$\{u\in\R^{n+1}:u^{\top}Wu<0\}$ contains a rational vector.  Scaling that
vector gives $v\in\Z^{n+1}$ with $v^{\top}Wv<0$.  Thus,
\[
 \bigcap_{v\in\Z^{n+1}}\{\eta:v^{\top}W(\eta)v\ge0\}=\mathcal P_n,
 \qquad
 \mathcal C_n\subseteq\mathcal P_n.
\]

We next introduce the second family of inequalities considered in this paper, namely the split inequalities, and define the corresponding split-inequality closure. For a primitive $v\in\Z^n$ and $s\in\Z$, define the \emph{primitive split
quadratic}
\begin{equation}\label{eq:split-polynomial}
 \widehat q_{v,s}(z)
 :=\frac12\bigl(v^{\top}z-s\bigr)
             \bigl(v^{\top}z-s-1\bigr).
\end{equation}
Its lattice-zero set is the primitive-split Delaunay polyhedron
\begin{equation}\label{eq:split-Delaunay}
 D_{v,s}:=Z(\widehat q_{v,s})
 =\{z\in\Z^n:v^{\top}z\in\{s,s+1\}\}.
\end{equation}
After an integer-affine change of variables, $D_{v,s}$ becomes
$\{0,1\}\times\Z^{n-1}$, and its convex hull becomes
$[0,1]\times\R^{n-1}$.

The \emph{split-inequality closure} is
\begin{equation}\label{eq:Sn}
 \mathcal S_n:=
 \bigcap_{\substack{v\in\Z^n\ \mathrm{primitive}\\s\in\Z}}
 \{W\in\mathcal H_n:L_{\widehat q_{v,s}}(W)\ge0\}.
\end{equation}
In matrix coordinates, the same inequality is
\begin{equation}\label{eq:split-matrix-form}
 v^{\top}Xv-(2s+1)v^{\top}x+s(s+1)\ge0,
\end{equation}
where we omit the positive factor $1/2$.  These inequalities are the
linearizations of the lattice disjunctions
$v^{\top}z\le s$ or $v^{\top}z\ge s+1$.  Burer and Letchford call
inequalities of this form \emph{split inequalities}
\cite[Section~5.2, Proposition~12]{BurerLetchford2014}.  We call the
intersection over primitive directions in \eqref{eq:Sn} the
\emph{split-inequality closure}.

\begin{proposition}[Basic containment relations]
\label{prop:splits-imply-psd}
For every $n\ge1$,
\begin{equation}\label{eq:universal-chain}
 \MIQ_n\subseteq\mathcal C_n\subseteq\mathcal S_n\subsetneq\mathcal P_n.
\end{equation}
\end{proposition}

\begin{proof}
Every lattice moment matrix belongs to $\mathcal C_n$. Indeed, for each
admissible $c$, the integer $c^{\top}\eta(z)$ is at least
$\beta_{\mathcal P_n}(c)$ and therefore at least
$\lceil\beta_{\mathcal P_n}(c)\rceil$. Because $\mathcal C_n$ is an
intersection of closed halfspaces, it is closed and convex. It therefore
contains the closed convex hull of the lattice generators, namely $\MIQ_n$.

Every split inequality is a rank-one CG cut. The primitive split quadratic
is integer-valued and has continuous minimum $-1/8$, as follows by
completing the square. Write
$\widehat q_{v,s}(z)=\alpha+c^{\top}\eta(z)$. By
Lemma~\ref{lem:integral-moment-coordinates}, $\alpha\in\Z$ and
$c\in\Z^{d_n}$. Moreover,
\[
 \beta_{\mathcal P_n}(c)
 =\min_{y\in\R^n}(\widehat q_{v,s}(y)-\alpha)
 =-\alpha-\frac18,
\]
so rounding gives $c^{\top}\eta\ge-\alpha$. This is exactly
$L_{\widehat q_{v,s}}\ge0$. Thus every split inequality occurs in the CG
family, giving $\mathcal C_n\subseteq\mathcal S_n$.

Finally, Burer and Letchford \cite[Theorem~4]{BurerLetchford2014} establish
$\mathcal S_n\subseteq\mathcal P_n$ using all integer split normals.
Primitive normals suffice because
for $W=\left[\begin{smallmatrix}1&x^{\top}\\x&X\end{smallmatrix}\right]
\in\mathcal S_n$, choosing $s=\lfloor v^{\top}x\rfloor$ in
\eqref{eq:split-matrix-form} gives $v^{\top}(X-xx^{\top})v\ge0$ for every
primitive $v\in\Z^n$. Every nonzero rational direction is a multiple of
such a vector, so continuity gives $X-xx^{\top}\succeq0$. The Schur
complement then gives $W\succeq0$. Strictness follows because
$W_{e_1/2}\in\mathcal P_n$, but
\eqref{eq:split-matrix-form} with $v=e_1$ and $s=0$ evaluates to
$1/4-1/2=-1/4<0$.
\end{proof}

We conclude by recording the integer-affine symmetry shared by all four sets.
Let
$T(z)=Uz+t$, where $U\in\GL_n(\Z)$ is unimodular (equivalently,
$U\in\Z^{n\times n}$ and $\det(U)=\pm1$) and $t\in\Z^n$, and define
\[
 R_T:=\begin{bmatrix}1&0\\t&U\end{bmatrix},
 \qquad
 \Phi_T(W):=R_TWR_T^{\top}.
\]
Then $\Phi_T(W_z)=W_{T(z)}$.

\begin{proposition}[Integer-affine symmetries]\label{prop:affine-symmetry}
The map $\Phi_T$ is an affine automorphism of $\mathcal H_n$ and satisfies
\[
 \Phi_T(\MIQ_n)=\MIQ_n,
 \qquad
 \Phi_T(\mathcal P_n)=\mathcal P_n,
 \qquad
 \Phi_T(\mathcal C_n)=\mathcal C_n,
 \qquad
 \Phi_T(\mathcal S_n)=\mathcal S_n.
\]
For every $p\in\mathcal Q_n$ and every $W\in\mathcal H_n$,
\begin{equation}\label{eq:affine-evaluation}
 L_{p\circ T^{-1}}\bigl(\Phi_T(W)\bigr)=L_p(W).
\end{equation}
Composition by $T^{-1}$ is an automorphism of
$\operatorname{Int}_2(\Z^n)$.
\end{proposition}

\begin{proof}
Because $T$ bijects $\Z^n$, the identity
$\Phi_T(W_z)=W_{T(z)}$ shows that $\Phi_T$ permutes the lattice moment
matrices.  Hence it preserves $\MIQ_n$.  The matrix $R_T$ is invertible, and
congruence by an invertible matrix preserves positive semidefiniteness.  Thus
$\Phi_T$ also preserves $\mathcal P_n$.

Identity \eqref{eq:affine-evaluation} holds at every $W_z$ by the definition
of linearization. Both sides are affine on $\mathcal H_n$, so the
affine-spanning consequence of Lemma~\ref{lem:lattice-zero-polynomial}
extends the identity to all of $\mathcal H_n$.

Since $T$ bijects $\Z^n$, a quadratic $p$ is integer-valued on $\Z^n$ if and
only if $p\circ T^{-1}$ is.  Since $T$ bijects $\R^n$, composition by
$T^{-1}$ preserves the infimum.

We now check the effect on a CG inequality.  For $c\in\Z^{d_n}$, write
\[
 c^{\top}\eta(T^{-1}(z))=\alpha+c'^{\top}\eta(z).
\]
Lemma~\ref{lem:integral-moment-coordinates} gives $\alpha\in\Z$ and
$c'\in\Z^{d_n}$.  Applying the same construction to $T$ gives the inverse
map, so $c\mapsto c'$ is a bijection of $\Z^{d_n}$.  Moreover,
\[
 \beta_{\mathcal P_n}(c')
 =\inf_{z\in\R^n}\bigl(c^{\top}\eta(T^{-1}(z))-\alpha\bigr)
 =\beta_{\mathcal P_n}(c)-\alpha.
\]
Thus $c$ has a finite relaxation value if and only if $c'$ does.  For these
normals, $\alpha$ is integral, so
\[
 \left\lceil\beta_{\mathcal P_n}(c')\right\rceil
 =\left\lceil\beta_{\mathcal P_n}(c)\right\rceil-\alpha.
\]
Together with \eqref{eq:affine-evaluation}, this implies that
$c^{\top}\eta(W)\ge\left\lceil\beta_{\mathcal P_n}(c)\right\rceil$
if and only if 
 $c'^{\top}\eta(\Phi_T(W))
\ge\left\lceil\beta_{\mathcal P_n}(c')\right\rceil$.
Thus $\Phi_T$ permutes the CG inequalities in \eqref{eq:Cn} and preserves
$\mathcal C_n$.

Finally,
\[
 \widehat q_{v,s}\circ T^{-1}
 =\widehat q_{U^{-\top}v,\,s+v^{\top}U^{-1}t}.
\]
The new normal is primitive because $U$ is unimodular, and the new level is
an integer.  Thus composition permutes the primitive split quadratics, so
$\Phi_T$ preserves $\mathcal S_n$.
\end{proof}

%% file: sections/03_facets.tex
\section{Facets of \texorpdfstring{$\MIQ_n$}{MIQ\_n}}
\label{sec:facets}

This section characterizes the facets of $\MIQ_n$ and determines when they
are bounded. The analysis uses Erdahl's structure theorem for
lattice-nonnegative quadratics to classify the facet-defining inequalities
and describe the geometry of unbounded perfect Delaunay polyhedra.

The affine-spanning consequence of Lemma~\ref{lem:lattice-zero-polynomial} shows that $\MIQ_n$ is full-dimensional in
$\mathcal H_n$, with dimension $d_n$. A \emph{facet} of $\MIQ_n$ is a nonempty
proper face of affine dimension $d_n-1$. Every facet is exposed. To see this,
let $F$ be a facet, choose $W\in\relint F$, and take a supporting hyperplane
of $\MIQ_n$ in $\mathcal H_n$ at $W$. The exposed face contains $F$ and has
dimension at most $d_n-1$, so it equals $F$. Thus every facet is the equality
set of a nonconstant valid affine inequality on $\mathcal H_n$.

\subsection{The facet--perfect-Delaunay correspondence}
\label{subsec:facet-correspondence}

The goal of this subsection is to show that the facets of $\mathrm{MIQ}_n$
are exactly the inequalities arising from perfect Delaunay polyhedra. We first
translate perfection into a dimension condition on the corresponding lifted
moment points. We then use Erdahl's classification of the extreme rays of
$\Erd(n)$ to prove the converse: every facet-defining inequality must arise
from a perfect Delaunay polyhedron.
\begin{lemma}
\label{lem:contact-rank}
For a Delaunay polyhedron $D\subseteq\Z^n$, the following statements are
equivalent:
\begin{enumerate}[label=\textup{(\roman*)},leftmargin=2.4em]
 \item $D$ is perfect, that is, $\dim\Space(D)=1$;
 \item
 \[
  \dim\lin\{W_z:z\in D\}=d_n;
 \]
 \item
 \[
  \dim\aff\{\eta(z):z\in D\}=d_n-1.
 \]
\end{enumerate}
\end{lemma}

\begin{proof}
Under the linear isomorphism $p\mapsto M_p$ between $\mathcal Q_n$ and
$\Smat^{n+1}$, 
\[
 \Space(D)
 =\{p\in\mathcal Q_n:\ip{M_p}{W_z}=0\text{ for every }z\in D\}
\]
and hence is the orthogonal complement of $\lin\{W_z:z\in D\}$.  Since
$\dim\mathcal Q_n=d_n+1$, rank--nullity shows that \textup{(i)} and
\textup{(ii)} are equivalent.

For every $z\in D$, the matrix $W_z$ satisfies $(W_z)_{00} = 1$. Thus, the matrices $W_z$, for $z\in D$, lie in an affine hyperplane that does not contain the origin. Their linear span has
dimension one more than their affine hull.  Finally, the affine bijection
$W\mapsto\eta(W)$ preserves affine dimension.  This proves the equivalence of
\textup{(ii)} and \textup{(iii)}.
\end{proof}

We use the following form of Erdahl's extreme-ray classification
\cite[Structure Theorem~2.1]{Erdahl1992}.
\begin{theorem}[Erdahl's extreme-ray theorem]
\label{thm:Erdahl-extreme-rays}
Every extreme ray of $\Erd(n)$ has a generator of exactly one of the
following types:
\begin{enumerate}[label=\textup{(\arabic*)},leftmargin=2.4em]
 \item the constant polynomial $1$;
 \item an affine square
 \[
  q(z)=(a^{\top}z+\beta)^2,
 \]
 where $a\in\R^n$ is not collinear with any nonzero vector of $\Z^n$;
 \item the defining quadratic $q_D$ of a perfect Delaunay polyhedron
 $D\subseteq\Z^n$.
\end{enumerate}
\end{theorem}

\begin{theorem}[Facet--perfect-Delaunay correspondence]
\label{thm:facet-correspondence}
The assignment
\[
 D\longmapsto
 F_D:=\MIQ_n\cap\{W\in\mathcal H_n:L_{q_D}(W)=0\}
\]
is a bijection between the perfect Delaunay polyhedra in $\Z^n$ and the
facets of $\MIQ_n$.  Equivalently, a nonconstant valid inequality
\[
 \ip{A}{X}+2b^{\top}x+\alpha\ge0
\]
defines a facet of $\MIQ_n$ if and only if the quadratic
$q(z)=z^{\top}Az+2b^{\top}z+\alpha$ is a positive multiple of the
defining quadratic of a perfect Delaunay polyhedron.
Moreover, for every integer-affine bijection $T$, the map $\Phi_T$ satisfies
$\Phi_T(F_D)=F_{T(D)}$.
\end{theorem}

\begin{proof}
Let $D$ be perfect.  Since $q_D\in\Erd(n)$, the inequality
$L_{q_D}\ge0$ is valid for $\MIQ_n$.  For every $z\in D$,
$L_{q_D}(W_z)=q_D(z)=0$, so $W_z\in F_D$.  By
Lemma~\ref{lem:contact-rank}, these contact matrices have affine hull of
dimension $d_n-1$.  Hence $\dim F_D\ge d_n-1$.  Moreover, $q_D\ne0$, so
$F_D$ lies in a proper supporting hyperplane of the $d_n$-dimensional set
$\MIQ_n$.  Therefore, $\dim F_D=d_n-1$, and $F_D$ is a facet.

Conversely, let $F$ be a facet and choose a nonconstant valid inequality
$L_q\ge0$ exposing it.  By \eqref{eq:moment-validity},
$q\in\Erd(n)$.  We first show that $q$ generates an extreme ray.  Suppose
$q=q_1+q_2$ with $q_1,q_2\in\Erd(n)$.  For every $W\in F$, the nonnegative
quantities $L_{q_1}(W)$ and $L_{q_2}(W)$ sum to $L_q(W)=0$.  Hence both
functionals vanish on $F$.  Since $\aff(F)$ is a hyperplane in the
$d_n$-dimensional affine space $\mathcal H_n$, every affine functional on
$\mathcal H_n$ that vanishes on $F$ is a scalar multiple of $L_q$.  The map
$p\mapsto L_p|_{\mathcal H_n}$ is injective by
Lemma~\ref{lem:lattice-zero-polynomial}, so $q_i=\lambda_iq$ for some
$\lambda_i\in\R$.  Because $F$ is a proper face, we may choose
$W^0\in\MIQ_n\setminus F$.  Then $L_q(W^0)>0$, and validity gives
$0\le L_{q_i}(W^0)=\lambda_iL_q(W^0)$.  Hence $\lambda_i\ge0$, so $q$
generates an extreme ray of $\Erd(n)$.

The first two types in Theorem~\ref{thm:Erdahl-extreme-rays} cannot define a
facet.  The constant inequality $1\ge0$ has an empty equality set.  For an
affine square $q(z)=(a^{\top}z+\beta)^2$, put
$v:=(\beta,a)\in\R^{n+1}$.  Every matrix in $\MIQ_n$ is positive
semidefinite, so
\[
 L_q(W)=v^{\top}Wv=0\quad\Longrightarrow\quad Wv=0.
\]
For $v\ne0$, consider the linear map
$T_v:\Smat^{n+1}\to\R^{n+1}$ given by $T_v(W)=Wv$.  This map is
surjective.  Indeed, for any $y\in\R^{n+1}$, the symmetric matrix
\[
 W_y:=\frac{yv^{\top}+vy^{\top}}{\|v\|_2^2}
      -\frac{v^{\top}y}{\|v\|_2^4}vv^{\top}
\]
satisfies $W_yv=y$.  Hence
$\dim\ker T_v=(d_n+1)-(n+1)=d_n-n$.  The equality face of $L_q\ge0$, if
nonempty, therefore has affine dimension at most $d_n-n$.  For $n\ge2$ this
is at most $d_n-2$, whereas a facet has dimension $d_n-1$.  When $n=1$,
every nonzero $a\in\R$ is collinear with the integer vector $1$, so the
affine-square type in Erdahl's theorem does not occur.  If $a=0$, the square
is constant.

Only the third type in Theorem~\ref{thm:Erdahl-extreme-rays} remains.  Thus
$q$ is a positive multiple of $q_D$ for some perfect Delaunay polyhedron
$D$, and $F=F_D$.

It remains to prove uniqueness.  Suppose $F_D=F_{D'}$.  The affine hull of a
facet is its unique supporting
hyperplane, so $L_{q_D}$ and $L_{q_{D'}}$ are proportional on $\mathcal H_n$.
The proportionality factor is positive because both inequalities are valid
and nonconstant.  Injectivity of the linearization map gives that $q_D$ and
$q_{D'}$ are positive scalar multiples, and therefore
\[
 D=Z(q_D)=Z(q_{D'})=D'.
\]

Finally, we prove equivariance.  Let $T(z)=Uz+t$ be an integer-affine
bijection.  Composition by
$T^{-1}$ maps $\Space(D)$ isomorphically onto $\Space(T(D))$, so perfection
is preserved.  Proposition~\ref{prop:affine-symmetry} gives
\[
 L_{q_D\circ T^{-1}}\bigl(\Phi_T(W)\bigr)=L_{q_D}(W),
\]
and $q_D\circ T^{-1}$ defines $T(D)$.  Consequently,
$\Phi_T(F_D)=F_{T(D)}$, proving equivariance.
\end{proof}

Given the individual facets classified by the last theorem, the next proposition shows
that their halfspaces describe the entire moment hull, even when the hull is
nonpolyhedral.

\begin{proposition}[Perfect-Delaunay facet inequalities describe the moment hull]
\label{prop:perfect-facets-describe-hull}
One has
\begin{equation}\label{eq:perfect-facet-description}
 \MIQ_n
 =\bigcap_{D\ \mathrm{perfect}}
   \{W\in\mathcal H_n:L_{q_D}(W)\ge0\}.
\end{equation}
\end{proposition}

\begin{proof}
Every perfect-Delaunay inequality is valid for $\MIQ_n$, so the inclusion
``$\subseteq$'' is immediate.  Let $K$ denote the intersection on the
right-hand side of \eqref{eq:perfect-facet-description}, and fix $W\in K$.

Every primitive-split Delaunay polyhedron is perfect
\cite{Erdahl1992}.  Therefore, the primitive-split
inequalities occur among the halfspaces defining $K$, and
\[
 K\subseteq\mathcal S_n\subseteq\mathcal P_n
\]
by Proposition~\ref{prop:splits-imply-psd}.

We next use the extreme rays to recover all inequalities in $\Erd(n)$.
The cone $\Erd(n)$ is closed, being an intersection of closed halfspaces. It is also pointed because if both $q$ and $-q$ belong to $\Erd(n)$, then $q$ vanishes on
$\Z^n$ and hence is the zero polynomial by
Lemma~\ref{lem:lattice-zero-polynomial}.  Every closed pointed cone in a
finite-dimensional space is generated by its extreme rays.

We now check the three types of extreme rays at $W$.  The constant inequality
is vacuous.  Every affine-square inequality is nonnegative at $W$ because
$W\in\mathcal P_n$.  Every perfect-Delaunay inequality holds by the
definition of $K$.  Thus $W$ satisfies the inequality associated with every
extreme ray in Theorem~\ref{thm:Erdahl-extreme-rays}.  Conic generation then
gives $L_q(W)\ge0$ for every $q\in\Erd(n)$.

Suppose, for contradiction, that $W\notin\MIQ_n$.  Since $\MIQ_n$ is closed
and convex, the separating-hyperplane theorem within $\mathcal H_n$ gives a
nonconstant affine inequality that is valid for $\MIQ_n$ and strictly
violated by $W$.  By \eqref{eq:moment-validity}, the corresponding
quadratic belongs to $\Erd(n)$.  This contradicts the preceding paragraph.
Therefore, $W\in\MIQ_n$, which proves the reverse inclusion.
\end{proof}

\subsection{Structure of perfect Delaunay polyhedra}
\label{subsec:bounded-facets}

The following coordinate form of Erdahl's structure theorem
\cite[Structure Theorem~2.1, p.~392]{Erdahl1992} describes a perfect Delaunay
polyhedron as a product of a bounded configuration and a free lattice.

\begin{theorem}[Erdahl's structure theorem]
\label{thm:perfect-delaunay-structure}
Let $D\subseteq\Z^n$ be a perfect Delaunay polyhedron. There exist an integer
$k\in\{1,\ldots,n\}$, an integer-affine bijection $T(z)=Uz+t$ with
$U\in\GL_n(\Z)$ and $t\in\Z^n$, and a finite perfect Delaunay polyhedron
$D_0\subseteq\Z^k$ with $\aff(D_0)=\R^k$, such that
\[
 T(D)=D_0\times\Z^{n-k}.
\]
Moreover, there is a quadratic $q_0\in\Erd(k)$ with $Z(q_0)=D_0$ such that
\[
 q_D\bigl(T^{-1}(u,v)\bigr)=q_0(u)
 =u^{\top}A_0u+2b_0^{\top}u+\alpha_0,
 \qquad A_0\succ0,
\]
for every $(u,v)\in\R^k\times\R^{n-k}$. In particular, $D$ is bounded if
and only if $k=n$.
\end{theorem}

This coordinate form follows by choosing bases of Erdahl's complementary
lattices and rescaling the defining quadratic of the bounded factor.

Since $A_0\succ0$, the quadratic $q_0$ has a finite, attained continuous
minimum. The same is therefore true of $q_D$, even when $D$ is unbounded.
This supplies the continuous minimum used to define depth in
Section~\ref{sec:depth}.

\begin{corollary}[Bounded facets]
\label{cor:bounded-facets}
For every perfect Delaunay polyhedron $D\subseteq\Z^n$, the corresponding
facet $F_D$ of $\MIQ_n$ is bounded if and only if $D$ is finite.
\end{corollary}

\begin{proof}
If $D$ is unbounded, the matrices $W_z$ with $z\in D$ form an unbounded
subset of $F_D$, so $F_D$ is unbounded. Suppose now that $D$ is finite.
Theorem~\ref{thm:perfect-delaunay-structure} implies that the quadratic
part of $q_D$ is positive definite. Completing the square, write
\[
 q_D(z)=(z-c)^{\top}A(z-c)-r,
 \qquad A\succ0.
\]
For $W=\left[\begin{smallmatrix}1&x^{\top}\\x&X\end{smallmatrix}\right]
\in F_D$, put $\Sigma=X-xx^{\top}\succeq0$, using
$\MIQ_n\subseteq\mathcal P_n$. Since $L_{q_D}(W)=0$, with
$\lambda=\lambda_{\min}(A)>0$ we obtain
\[
 r=(x-c)^{\top}A(x-c)+\ip{A}{\Sigma}
 \ge\lambda\bigl(\|x-c\|_2^2+\operatorname{tr}\Sigma\bigr).
\]
Thus $x$ and $\operatorname{tr}\Sigma$ are bounded. Positive
semidefiniteness bounds every entry of $\Sigma$, so $X=xx^{\top}+\Sigma$
is also bounded. Hence $F_D$ is bounded.
\end{proof}

Thus the bounded facets correspond exactly to perfect Delaunay polytopes.
We next determine which perfect-Delaunay facet inequalities are recovered
by one round of CG cuts.

%% file: sections/04_depth.tex
\section{Depth and the rank-one Chv\'atal--Gomory closure}
\label{sec:depth}

In Section~\ref{sec:facets}, we characterized the facet-defining inequalities of $\MIQ_n$.  We now explore which of these inequalities can be obtained as rank-one CG cuts from $\mathcal P_n$.  Given that CG rounding depends on positive scaling, we first normalize the defining quadratics so that their linearizations have primitive integer normals in integral moment coordinates.  We then introduce the notion of depth for quadratics and perfect Delaunay polyhedra, with which we categorize perfect Delaunay polyhedra as `shallow' or `deep'.  We show that, when a perfect Delaunay polyhedron is shallow, its associated facet-defining inequality for $\MIQ_n$ can be obtained in the first round of CG cuts for $\mathcal P_n$. On the other hand, we prove that the facet-defining inequalities of $\MIQ_n$ associated with deep perfect Delaunay polyhedra cannot be recovered from the rank-one CG closure of $\mathcal P_n$. 

Determining whether CG rounding from $\mathcal P_n$ recovers a facet-defining inequality requires an integral moment normal on its defining ray. For a nonconstant quadratic $q\in\mathcal Q_n$, let $\R_+q:=\{tq:t\ge0\}$ denote the ray generated by $q$.  We call this ray \emph{rational} if it contains a nonzero quadratic with integer Newton coefficients.

\begin{lemma}
\label{lem:perfect-ray-rationality}
For every perfect Delaunay polyhedron $D\subseteq\Z^n$, the ray $\R_+q_D$ is rational. 
\end{lemma}

\begin{proof}
By perfection, the system
\[
 \alpha+c^{\top}\eta(z)=0
\]
for $z\in D$ has a one-dimensional solution space, consisting of the Newton coefficient
vectors of the polynomials in $\Space(D)=\R q_D$.  Since there are only
$d_n+1$ unknowns, finitely many of these equations define the same solution
space.  Their coefficients are integers because $\eta(z)$ is integral for
every $z\in D$.  Gaussian elimination over $\Q$ therefore gives a nonzero
rational solution.  Choosing its sign and clearing denominators yields a
positive multiple of $q_D$ with integer Newton coefficients.
\end{proof}

Rationality allows us to choose an integral representative of the defining ray. We next normalize this representative by requiring its nonconstant Newton coefficient vector to be primitive. A nonzero coefficient vector $c\in\Z^{d_n}$ is \emph{primitive} if its entries have greatest common divisor (gcd) one. The next lemma shows that this normalization is unique and that the constant coefficient is also integral when the quadratic has a lattice zero.

\begin{lemma}
\label{lem:primitive-normalization}
Suppose that the nonconstant ray $\R_+q$ is rational and that $q$ has a lattice zero.  Then there is a unique positive scaling
\[
 \widehat q(z)=\alpha+c^{\top}\eta(z),
 \qquad
 \alpha\in\Z,\qquad c\in\Z^{d_n}\text{ primitive}.
\]
\end{lemma}

\begin{proof}
Since the ray is rational, choose a positive multiple
\[
p(z)=\alpha'+{c'}^\top\eta(z)
\]
of \(q\) with \(\alpha'\in\mathbb Z\) and \(c'\in\mathbb Z^{d_n}\). Since \(q\) is nonconstant, \(c'\neq 0\). Let
\[
g:=\gcd(c'_1,\ldots,c'_{d_n})
\]
and define
\[
\widehat q:=\frac{p}{g}
=\alpha+c^\top\eta(z),
\qquad
c:=\frac{c'}{g}.
\]
By construction, \(c\in\mathbb Z^{d_n}\) is primitive.

It remains to show that \(\alpha=\alpha'/g\) is an integer. Let \(z^0\in\mathbb Z^n\) be a lattice zero of \(q\). Since \(p\) is a positive multiple of \(q\), \(p(z^0)=0\), and therefore
\[
\alpha'
=-{c'}^\top\eta(z^0).
\]
Because every entry of \(c'\) is divisible by \(g\), the right-hand side is divisible by \(g\). Hence \(\alpha'/g\in\mathbb Z\), so \(\widehat q\) is integer-valued by Lemma~\ref{lem:integral-moment-coordinates}.

Finally, this scaling is unique. Indeed, any other positive scaling of \(q\) with a primitive integral nonconstant Newton coefficient vector must have its coefficient vector equal to a positive scalar multiple of \(c\). Since both vectors are integral and primitive, that scalar must be one.
\end{proof}

We call $\widehat q$ the \emph{primitive integer-valued generator} of the ray.  If $q$ is nonconstant, integer-valued, and has a lattice zero, let $\Delta(q)$ denote the greatest common divisor of the entries of $c_q$.  Since $q$ itself is an integral representative of its defining ray, Lemma~\ref{lem:primitive-normalization} gives
\begin{equation}\label{eq:primitive-content}
 \widehat q=\frac{q}{\Delta(q)}.
\end{equation}

For a perfect Delaunay polyhedron $D$, the normalized inequality $L_{\widehat q_D}(W)\ge0$ defines the same facet of $\MIQ_n$ as $L_{q_D}(W)\ge0$.  The minimum value of $L_{\widehat q_D}$ over $\MIQ_n$ is zero, while by \eqref{eq:beta-formula}, its minimum over $\mathcal P_n$ is the continuous minimum of $\widehat q_D$. The notion of depth introduced below quantifies the gap between these two minima.

\begin{definition}
\label{def:depth}
Let $q$ be a nonconstant lattice-nonnegative quadratic such that $ \R_+q$ is rational, $q$ has at least one lattice-zero, and  $\min_{y\in\R^n}q(y)>-\infty$. The \emph{depth} of $q$ is
\begin{equation}\label{eq:depth}
 \delta(q):=-\min_{y\in\R^n}\widehat q(y)
 =\min_{z\in\Z^n}\widehat q(z)
  -\min_{y\in\R^n}\widehat q(y)\ge0.
\end{equation}
This quantity depends only on the positive ray.  For a perfect Delaunay polyhedron $D$, we write $\delta(D):=\delta(q_D)$ and call $D$ \emph{shallow} when $\delta(D)<1$ and \emph{deep} when $\delta(D)\ge1$.
\end{definition}

The lattice minimum in \eqref{eq:depth} is zero because $q$ is lattice-nonnegative and has a lattice zero.  For a perfect Delaunay polyhedron $D$, Lemma~\ref{lem:perfect-ray-rationality} ensures that its defining ray is rational, while the representation in Theorem~\ref{thm:perfect-delaunay-structure} ensures that the continuous minimum of $q_D$ is finite and attained. Hence $\delta(D)$ is well defined.

Equation~\eqref{eq:beta-formula} gives
\begin{equation}\label{eq:depth-relaxation-bounds}
 \min_{W\in\MIQ_n}L_{\widehat q_D}(W)=0,
 \qquad
 \min_{W\in\mathcal P_n}L_{\widehat q_D}(W)=-\delta(D).
\end{equation}
Thus, the depth of $D$ is the difference between the minimum values of $L_{\widehat q_D}$ over $\MIQ_n$ and $\mathcal P_n$.  If $q$ in Definition~\ref{def:depth} is integer-valued, then \eqref{eq:primitive-content} gives
\[
 \delta(q)=\frac{-\min_{y\in\R^n}q(y)}{\Delta(q)}.
\]

The following example shows how normalization affects CG rounding. \\

\noindent\textbf{Example.} 
The seven-dimensional Gosset polytope $3_{21}$ is a perfect Delaunay
polytope \cite{DutourErdahlRybnikov2007}.  We use the following defining
quadratic on $\R^7$, written in the Newton basis:
\[
 \begin{aligned}
 q_{G_7}(z)={}&
 8\sum_{i=1}^{5}\binom{z_i}{2}
 +16\binom{z_6}{2}+12\binom{z_7}{2}+2z_6\\
 &+6\sum_{1\le i<j\le5}z_iz_j
 +(10z_6+8z_7)\sum_{i=1}^{5}z_i+12z_6z_7.
 \end{aligned}
\]
This polynomial has zero constant term.  Its nonconstant Newton
coefficients are all even, and the coefficient of $z_6$ is two.
Their greatest common divisor is therefore two, so
$\widehat q_{G_7}=q_{G_7}/2$ is the primitive integer-valued generator.
The continuous minimum of $q_{G_7}$ is $-3/2$, attained at $e_7/2$, where
$e_7$ is the seventh unit vector.  Appendix~\ref{app:low-dimensional}
gives the matrix representation and positive-definiteness check.
Thus CG rounding gives
\[
 \begin{aligned}
 L_{q_{G_7}}(W)\ge-\frac32
 &\quad\xrightarrow{\mathrm{CG}}\quad L_{q_{G_7}}(W)\ge-1,\\
 L_{\widehat q_{G_7}}(W)\ge-\frac34
 &\quad\xrightarrow{\mathrm{CG}}\quad L_{\widehat q_{G_7}}(W)\ge0
 \end{aligned}.
\]
The two polynomials define the same facet, but rounding with the primitive
normal recovers its exact inequality.  Using $q_{G_7}$ directly gives only
the weaker bound $L_{\widehat q_{G_7}}(W)\ge-1/2$.
The correctly normalized depth is $\delta(3_{21})=3/4$.

To use depth for all perfect Delaunay polyhedra, we next examine how it behaves under the structural operations introduced in Section~\ref{subsec:bounded-facets}.  The following proposition shows that depth is positive and invariant under integer-affine transformations and products with free lattice factors.

\begin{proposition}
\label{prop:depth-properties}
Every perfect Delaunay polyhedron $D$ has strictly positive depth. In addition, if $T$ is an integer-affine bijection of $\R^n$ and $k\in\Z_{++}$, then $T(D)$ and $D\times\Z^k$ are perfect. Moreover,
\[
 \delta(T(D))=\delta(D)
 \qquad\text{and}\qquad
 \delta(D\times\Z^k)=\delta(D).
\]
\end{proposition}

\begin{proof}
We first prove $\delta(D)>0$ by showing that $\delta(D)=0$ would force $\widehat q_D\equiv0$.  Suppose that $\delta(D)=0$, so $\min_{\R^n}\widehat q_D=0$.  Write $\widehat q_D(y)=y^{\top}Ay+2b^{\top}y+\alpha$, with symmetric $A$, and choose $z^0\in D$.  The equality $\widehat q_D(z^0)=0$ makes $z^0$ a global minimizer.  Thus, $A\succeq0$ and $Az^0+b=0$, which give
\[
 \widehat q_D(y)=(y-z^0)^{\top}A(y-z^0).
\]
For a positive-semidefinite matrix $A$, the right-hand side is zero exactly when $y-z^0\in\ker A$.  Therefore,
\[
 \{y\in\R^n:\widehat q_D(y)=0\}=z^0+\ker A.
\]
Every $z\in D$ satisfies $\widehat q_D(z)=0$, so $D\subseteq z^0+\ker A$.  The affine subspace $z^0+\ker A$ must therefore contain $\aff(D)=\R^n$, as established in Section~\ref{sec:delaunay}.  It follows that $\ker A=\R^n$ and $A=0$.  The displayed expression for $\widehat q_D$ then gives $\widehat q_D\equiv0$, contradicting that it is a nonzero generator of the defining ray. 

Next, we show that $T(D)$ is perfect and that $\delta(T(D))=\delta(D)$.
By Theorem~\ref{thm:facet-correspondence}, $T(D)$ is perfect, and its
defining ray is generated by $\widehat q_D\circ T^{-1}$.  This polynomial
is integer-valued and vanishes at the integer point $T(z^0)$, where
$z^0\in D$.

To compute the depth of $T(D)$, we must check that this generator's
nonconstant Newton coefficient vector is primitive.  Suppose it is not,
so some integer $m>1$ divides all these coefficients.  Since the
polynomial vanishes at $T(z^0)$, its constant coefficient is also
divisible by $m$.  Thus, $(\widehat q_D\circ T^{-1})/m$ is integer-valued.
Composing with $T$ shows that $\widehat q_D/m$ is integer-valued.
Lemma~\ref{lem:integral-moment-coordinates} then implies that every entry
of its nonconstant Newton coefficient vector $c_{\widehat q_D}/m$ is an
integer.
This means that $m$ divides every entry of $c_{\widehat q_D}$,
contradicting that this vector has gcd one.  Hence,
$\widehat q_D\circ T^{-1}$ is the primitive integer-valued generator of
the transformed ray.  Because $T$ is a bijection of $\R^n$, the continuous
minimum is unchanged.  Therefore, $\delta(T(D))=\delta(D)$.

Finally, we argue that $D\times\Z^k$ is perfect and that $\delta(D\times\Z^k)=\delta(D)$.  The extension $(z,u)\mapsto\widehat q_D(z)$ has lattice-zero set $D\times\Z^k$.  For $p\in\Space(D\times\Z^k)$, write
\[
 p(z,u)=r(z)+u^{\top}Bu+\ell(z)^{\top}u,
 \qquad r\in\mathcal Q_n,\quad B\in\Smat^k,
\]
where $\ell:\R^n\to\R^k$ is affine.  For each $z\in D$, Lemma~\ref{lem:lattice-zero-polynomial} gives $p(z,u)=0$ for all $u\in\R^k$, so $B=0$, $\ell(z)=0$, and $r(z)=0$.  The identity $\aff(D)=\R^n$ gives $\ell\equiv0$, while $r\in\Space(D)=\R\widehat q_D$.  Thus, $\Space(D\times\Z^k)$ consists exactly of the scalar multiples of the extension, which indicates perfection.  Its nonconstant Newton coefficient vector is obtained by inserting zero entries into that of $\widehat q_D$, so it remains primitive.  Its continuous minimum is unchanged, so $\delta(D\times\Z^k)=\delta(D)$.
\end{proof}

We now determine when CG rounding strengthens the bound $-\delta(D)$ in \eqref{eq:depth-relaxation-bounds} to the exact bound zero, thereby giving the facet-defining inequality $L_{\widehat q_D}(W)\ge0$ of $\MIQ_n$.  We express rank-one CG cuts for $\mathcal P_n$ in the form $L_p(W)\ge0$, where $p$ is an integer-valued quadratic.  For a nonzero $c\in\Z^{d_n}$ with $\beta_{\mathcal P_n}(c)>-\infty$, we call 
\begin{equation}\label{eq:CG-cut-polynomial}
 p(z):=c^{\top}\eta(z)
       -\left\lceil\beta_{\mathcal P_n}(c)\right\rceil
\end{equation}
the \emph{rank-one CG cut polynomial associated with $c$}.  The associated inequality $L_p(W)\ge0$ is the rank-one CG cut $c^{\top}\eta(W)\ge\lceil\beta_{\mathcal P_n}(c)\rceil$ for $\mathcal P_n$.  Equation~\eqref{eq:beta-formula} gives
\begin{equation}\label{eq:CG-cut-minimum}
 \min_{y\in\R^n}p(y)
 =\beta_{\mathcal P_n}(c)
  -\left\lceil\beta_{\mathcal P_n}(c)\right\rceil\in(-1,0].
\end{equation}

In the next lemma, we discuss exactly when a nonzero integer-valued quadratic is a rank-one CG cut polynomial.

\begin{lemma}
\label{lem:CG-polynomial-dictionary}
A nonzero polynomial $p\in\operatorname{Int}_2(\Z^n)$ is a rank-one CG cut if and only if it has a finite continuous minimum satisfying
\begin{equation}\label{eq:intrinsic-CG-polynomials}
 -1<\min_{y\in\R^n}p(y)\le0.
\end{equation}
Its normal is uniquely determined by its nonconstant Newton coefficients. Thus,
\[
 \mathcal C_n
 =\bigcap_{\substack{
 p\in\operatorname{Int}_2(\Z^n)\setminus\{0\}\\
 -1<\min_{y\in\R^n}p(y)\le0}}
 \{W\in\mathcal H_n:L_p(W)\ge0\}.
\]
\end{lemma}

\begin{proof}
Equation~\eqref{eq:CG-cut-minimum} proves necessity.  For the converse, consider $p(z)=\alpha+c^{\top}\eta(z)$, where $\alpha\in\Z$ and $c\in\Z^{d_n}$.  A nonzero integer-valued constant cannot satisfy \eqref{eq:intrinsic-CG-polynomials}, so $c\ne0$.  If $m:=\min_{y\in\R^n}p(y)\in(-1,0]$, then \eqref{eq:beta-formula} gives
\[
 \beta_{\mathcal P_n}(c)=m-\alpha,
 \qquad
 \left\lceil\beta_{\mathcal P_n}(c)\right\rceil=-\alpha.
\]
Hence, $p$ has the form \eqref{eq:CG-cut-polynomial}. The Newton-basis representation uniquely determines the coefficient vector $c$.  We substitute this polynomial representation of the cuts into \eqref{eq:Cn} and apply the coordinate identification of Lemma~\ref{lem:integral-moment-coordinates} to obtain the asserted description of $\mathcal C_n$.  We omit $c=0$ because it gives only the trivial inequality.
\end{proof}

We next determine when a single CG cut along a facet normal recovers the
facet inequality.

\begin{proposition}[CG rounding along a facet normal]
\label{prop:facet-normal-rounding}
Let $D\subseteq\Z^n$ be a perfect Delaunay polyhedron, with primitive
integer-valued generator
\[
 \widehat q_D(z)=\alpha_D+c_D^{\top}\eta(z),
 \qquad \alpha_D\in\Z,\qquad c_D\in\Z^{d_n}\text{ primitive}.
\]
For every integer $\mu\ge1$, the rank-one CG cut obtained from the normal
$\mu c_D$, after division by $\mu$, is
\begin{equation}\label{eq:multiple-normal-cut}
 L_{\widehat q_D}(W)\ge
 -\frac{\lfloor\mu\delta(D)\rfloor}{\mu}.
\end{equation}
This cut is the facet inequality $L_{\widehat q_D}(W)\ge0$ if and only if
$\mu\delta(D)<1$.  Equivalently, $\mu\widehat q_D$ is a rank-one CG cut polynomial if and only if $\mu\delta(D)<1$.
\end{proposition}

\begin{proof}
Since
\[
\widehat q_D(z)=\alpha_D+c_D^\top\eta(z)
\qquad\text{and}\qquad
\min_{z\in\R^n}\widehat q_D(z)=-\delta(D),
\]
we have
\[
\beta_{\mathcal P_n}(\mu c_D)
=-\mu\alpha_D-\mu\delta(D).
\]
Because \(\mu\alpha_D\) is integral,
\[
\left\lceil\beta_{\mathcal P_n}(\mu c_D)\right\rceil
=-\mu\alpha_D-\lfloor\mu\delta(D)\rfloor.
\]
Hence the CG cut associated with the normal \(\mu c_D\) is
\[
\mu c_D^\top\eta(W)
\ge
-\mu\alpha_D-\lfloor\mu\delta(D)\rfloor.
\]
Adding \(\mu\alpha_D\) and dividing by \(\mu\) gives
\[
L_{\widehat q_D}(W)
\ge
-\frac{\lfloor\mu\delta(D)\rfloor}{\mu}.
\]
This is the facet inequality \(L_{\widehat q_D}(W)\ge0\) exactly when
\(\lfloor\mu\delta(D)\rfloor=0\), or equivalently when
\(\mu\delta(D)<1\).

Finally, the corresponding cut polynomial is
\[
\mu\widehat q_D+\lfloor\mu\delta(D)\rfloor,
\]
which equals \(\mu\widehat q_D\) exactly under the same condition.
\end{proof}

Taking $\mu=1$ shows that every shallow facet is obtained by CG rounding in
its primitive normal.  Because $c_D$ is primitive, every integral normal
on its positive ray has the form $\mu c_D$ with an integer $\mu\ge1$.
For a deep perfect Delaunay polyhedron $D$, Proposition~\ref{prop:facet-normal-rounding}
therefore rules out recovering the facet inequality by a single cut in this
normal direction.  It remains to determine whether rank-one CG cuts with other normals jointly imply this inequality.  In general, an infinite family of inequalities can imply an inequality without containing a positive multiple of it.  For example,
\[
 \bigcap_{k\ge1}\{(x,y)\in\R^2:kx+y\ge0\}
 \subseteq\{(x,y)\in\R^2:x\ge0\},
\]
because any $x<0$ violates $kx+y\ge0$ for sufficiently large $k$.
For the rigidity argument, we first present a coefficient bound for arbitrary
real quadratics bounded below by $-1$.

\begin{lemma}[Newton coefficient bound]
\label{lem:newton-coefficient-bound}
Let $p(z)=\alpha_p+c_p^{\top}\eta(z)\in\mathcal Q_n$ satisfy
$p(y)\ge-1$ for every $y\in\R^n$.  For every $W^0\in\mathcal P_n$ with
$W^0\succ0$,
\[
 \|c_p\|_\infty\le
 \frac{3\bigl(1+L_p(W^0)\bigr)}{\lambda_{\min}(W^0)}.
\]
\end{lemma}

\begin{proof}
The idea is to represent the globally nonnegative quadratic \(p+1\) by a positive-semidefinite matrix and then use \(W^0\succ0\) to control its entries. Let $M:=M_{p+1}$, the coefficient matrix of the shifted polynomial
$p+1$.  For $t\ne0$ and $y\in\R^n$,
\[
\begin{bmatrix}t\\y\end{bmatrix}^{\!\top}
M
\begin{bmatrix}t\\y\end{bmatrix}
=
t^2\bigl(p(y/t)+1\bigr)\ge0.
\]
By continuity, the same inequality holds when \(t=0\), and hence \(M\succeq0\). Let \(\lambda:=\lambda_{\min}(W^0)>0\). Since \(W^0\succeq\lambda I\), \(W^0_{00}=1\), and \(M\succeq0\),
\[
1+L_p(W^0)
=
\langle M,W^0\rangle
\ge
\lambda\operatorname{tr}(M).
\]
Therefore,
\[
\operatorname{tr}(M)
\le
\frac{1+L_p(W^0)}{\lambda}.
\]

It remains to relate the entries of \(M\) to the Newton coefficients of \(p\). Since adding a constant does not change the nonconstant Newton coefficients, \(p\) and \(p+1\) have the same vector \(c_p\). Solving the relations in \eqref{eq:newton-to-matrix} for these coefficients gives
\[
a_i=2M_{0i}+M_{ii},\qquad
d_i=2M_{ii},\qquad
\gamma_{ij}=2M_{ij}\quad(i<j).
\]
Because \(M\succeq0\), every entry of \(M\) has absolute value at most \(\operatorname{tr}(M)\). Hence
\[
|d_i|\le2\operatorname{tr}(M),\qquad
|\gamma_{ij}|\le2\operatorname{tr}(M),\qquad
|a_i|\le3\operatorname{tr}(M).
\]
Thus,
\[
\|c_p\|_\infty
\le
3\operatorname{tr}(M)
\le
\frac{3\bigl(1+L_p(W^0)\bigr)}
{\lambda_{\min}(W^0)}.
\]
\end{proof}

We use this bound to prove the following rigidity theorem.  Applied to
rank-one CG cuts, it rules out recovering a deep facet inequality through
cuts with other normals.

\begin{theorem}[Rigidity theorem]
\label{thm:rigidity-depth}
Let $D\subseteq\Z^n$ be a perfect Delaunay polyhedron, and let $\mathcal F\subseteq\operatorname{Int}_2(\Z^n)\setminus\{0\}$ satisfy $\inf_{\R^n}p>-1$ for every $p\in\mathcal F$.  Define
\[
 \mathcal K_{\mathcal F}
 :=\{W\in\mathcal P_n:L_p(W)\ge0\text{ for every }p\in\mathcal F\}.
\]
Then $L_{\widehat q_D}\ge0$ is valid for $\mathcal K_{\mathcal F}$ if and only if $\mathcal F$ contains $\mu\widehat q_D$ for some integer $\mu\ge1$. 
\end{theorem}

\begin{proof}
We first prove validity when $\mathcal F$ contains $\mu\widehat q_D$ for some integer $\mu\ge1$.  Every $p\in\mathcal F$ is lattice-nonnegative because its values on $\Z^n$ are integers strictly greater than $-1$.  Thus, $\MIQ_n\subseteq\mathcal K_{\mathcal F}$.  If $\mu\widehat q_D\in\mathcal F$, every $W\in\mathcal K_{\mathcal F}$ satisfies $\mu L_{\widehat q_D}(W)\ge0$.  We divide by $\mu>0$ and conclude that $L_{\widehat q_D}\ge0$ is valid for $\mathcal K_{\mathcal F}$. 

For the converse, we assume that $L_{\widehat q_D}\ge0$ is valid for $\mathcal K_{\mathcal F}$ and prove that $\mathcal F$ contains a positive integral multiple of $\widehat q_D$.  For a contradiction, suppose that $\mathcal F$ contains no such multiple.  We construct a matrix in $\mathcal K_{\mathcal F}$ at which $L_{\widehat q_D}$ is negative.  Specifically, we will first choose a positive-definite matrix on the facet associated with $D$ and show that it strictly satisfies every inequality indexed by $\mathcal F$.  We will then use a uniform coefficient estimate to choose a perturbation that remains in $\mathcal K_{\mathcal F}$ and makes $L_{\widehat q_D}$ negative. 

To construct the starting matrix on the facet associated with $D$, we take the average of the matrices $W_z$ corresponding to selected points $z\in D$.  By Lemma~\ref{lem:contact-rank}, we choose $z^1,\ldots,z^{d_n}\in D$ such that $W_{z^1},\ldots,W_{z^{d_n}}$ form a basis of $\operatorname{span}\{W_z:z\in D\}$, and define
\[
 \overline W:=\frac1{d_n}\sum_{i=1}^{d_n}W_{z^i}.
\]
For every $p\in\mathcal Q_n$, we have $p(z)=\ip{M_p}{W_z}$.  The basis property therefore gives
\[
 \{p\in\mathcal Q_n:p(z^i)=0\text{ for }i=1,\ldots,d_n\}
 =\Space(D)=\R\widehat q_D.
\]
As a convex combination of the matrices $W_{z^i}$, the matrix $\overline W$ belongs to $\MIQ_n$.  Moreover, $L_{\widehat q_D}(\overline W)=0$ because $\widehat q_D(z^i)=0$ for every selected point.  Thus, $\overline W$ belongs to the facet associated with $D$.

We now prove that $\overline W\succ0$.  For a contradiction, suppose that $\overline W$ is not positive definite.  Each $W_{z^i}$ is positive semidefinite, so $\overline W\succeq0$.  We can therefore choose a nonzero $v\in\R^{n+1}$ with
\[
 0=v^{\top}\overline Wv
  =\frac1{d_n}\sum_{i=1}^{d_n}
     \left(v^{\top}\begin{bmatrix}1\\z^i\end{bmatrix}\right)^2.
\]
Every summand is nonnegative, so $v^{\top}\left[\begin{smallmatrix}1\\z^i\end{smallmatrix}\right]=0$ for $i=1,\ldots,d_n$.  The preceding characterization of $\Space(D)$ then places the nonzero affine polynomial $z\mapsto v^{\top}\left[\begin{smallmatrix}1\\z\end{smallmatrix}\right]$ in $\R\widehat q_D$.  Because $\widehat q_D$ is not affine, as established in Section~\ref{sec:delaunay}, we reach a contradiction.  Hence, $\overline W\succ0$. 

Next, we establish a positive lower bound on $L_p(\overline W)$ that holds for every $p\in\mathcal F$.  Let
\[
 s(p):=L_p(\overline W)
      =\frac1{d_n}\sum_{i=1}^{d_n}p(z^i).
\]
The summands are nonnegative integers.  Fix $p\in\mathcal F$ and, for a contradiction, suppose that $s(p)=0$.  Then $p(z^i)=0$ for every $i$, so $p=\mu\widehat q_D$ for some nonzero real $\mu$.  Because $\widehat q_D$ is nonzero and lattice-nonnegative, Lemma~\ref{lem:lattice-zero-polynomial} gives a point $z\in\Z^n$ with $\widehat q_D(z)>0$.  The inequality $p(z)=\mu\widehat q_D(z)\ge0$ and the condition $\mu\ne0$ imply $\mu>0$.  Moreover, $c_p=\mu c_D$ is integral and $c_D$ is primitive, so we obtain $\mu\in\Z$ from B\'ezout's identity.  We have found a positive integral multiple of $\widehat q_D$ in $\mathcal F$, which is a contradiction.  Thus, at least one summand in the definition of $s(p)$ is a positive integer, and we obtain
\begin{equation}\label{eq:discrete-slack}
 s(p)\ge\frac1{d_n},\quad \forall\,p\in\mathcal F.
\end{equation}

To control how the inequalities change under perturbation, apply
Lemma~\ref{lem:newton-coefficient-bound} to each $p\in\mathcal F$ at
$W^0=\overline W$.  With $\lambda:=\lambda_{\min}(\overline W)>0$, this gives
\[
 \|c_p\|_\infty\le\frac3\lambda\bigl(1+s(p)\bigr).
\]
Since \eqref{eq:discrete-slack} implies
$1+s(p)\le(d_n+1)s(p)$, we obtain
\begin{equation}\label{eq:slope-slack}
 \|c_p\|_\infty\le C_Ds(p),\, \forall\, p\in\mathcal F,
 \qquad C_D:=\frac{3(d_n+1)}{\lambda}.
\end{equation}
The constant $C_D$ depends only on the selected points of $D$. 

Finally, we construct a matrix in $\mathcal K_{\mathcal F}$ that violates $L_{\widehat q_D}\ge0$.  We choose $\varepsilon>0$ sufficiently small that
\[
 W^\varepsilon:=W\bigl(\eta(\overline W)-\varepsilon c_D\bigr)\succ0,
 \qquad \varepsilon C_D\|c_D\|_1\le\frac12.
\]
We can make this choice because $\overline W\succ0$ and $\eta\mapsto W(\eta)$ is continuous.  We use \eqref{eq:evaluation-dictionary} and \eqref{eq:slope-slack} to verify that, for every $p\in\mathcal F$,
\[
 \begin{aligned}
 L_p(W^\varepsilon)
 &=s(p)-\varepsilon c_p^{\top}c_D\\
 &\ge s(p)-\varepsilon\|c_p\|_\infty\|c_D\|_1\\
 &\ge\bigl(1-\varepsilon C_D\|c_D\|_1\bigr)s(p)
 \ge\frac12s(p)>0.
 \end{aligned}
\]
Thus, $W^\varepsilon\in\mathcal K_{\mathcal F}$, whereas
\[
 L_{\widehat q_D}(W^\varepsilon)
 =-\varepsilon\|c_D\|_2^2<0.
\]
We assumed that $L_{\widehat q_D}\ge0$ is valid for $\mathcal K_{\mathcal F}$, so we reach a contradiction. 
\end{proof}

We now apply Theorem~\ref{thm:rigidity-depth} to the family of all rank-one CG cut polynomials to characterize the facet-defining inequalities that are valid for $\mathcal C_n$. 

\begin{corollary}
\label{cor:depth-criterion}
For every perfect Delaunay polyhedron $D\subseteq\Z^n$, the following conditions are equivalent:
\begin{enumerate}
\item[(i)] $\widehat q_D$ is a rank-one CG cut polynomial; 
\item[(ii)] $L_{\widehat q_D}\geq 0$ is valid for $\mathcal C_n$; 
\item[(iii)] $\delta(D)<1$. 
\end{enumerate}
\end{corollary}

The CG rank of an inequality relative to $\mathcal P_n$ is the minimum
number of successive CG rounds needed to recover it, or $+\infty$ if no
finite number suffices. For every deep $D$, the facet-defining inequality
$L_{\widehat q_D}(W)\ge0$ fails at some point of $\mathcal C_n$. Therefore,
its CG rank is at least two.

\begin{proof}
With $\mu=1$, Proposition~\ref{prop:facet-normal-rounding} gives the equivalence between the first and third conditions.  To establish the equivalence with validity for $\mathcal C_n$, let $\mathcal F$ be the family of all rank-one CG cut polynomials.  Their continuous minima belong to $(-1,0]$, so they satisfy the conditions of Theorem~\ref{thm:rigidity-depth}.  Lemma~\ref{lem:CG-polynomial-dictionary} and the inclusion $\mathcal C_n\subseteq\mathcal P_n$ from Section~\ref{sec:rank-one-closure} give $\mathcal K_{\mathcal F}=\mathcal C_n$.  By Theorem~\ref{thm:rigidity-depth}, the inequality $L_{\widehat q_D}\ge0$ is valid for $\mathcal C_n$ if and only if $\mathcal F$ contains $\mu\widehat q_D$ for some integer $\mu\ge1$.  By Proposition~\ref{prop:facet-normal-rounding}, such a multiple exists if and only if $\delta(D)<1$.
\end{proof}

We combine Corollary~\ref{cor:depth-criterion} with Proposition~\ref{prop:perfect-facets-describe-hull} to characterize exactly when $\mathcal C_n=\MIQ_n$. 

\begin{corollary}
\label{cor:depth-exact-recovery}
For every $n\ge1$, $\mathcal C_n=\MIQ_n$ if and only if $\delta(D)<1$ for every perfect Delaunay polyhedron $D\subseteq\Z^n$.
\end{corollary}

\begin{proof}
We first assume that every perfect Delaunay polyhedron in $\Z^n$ is shallow.  By Corollary~\ref{cor:depth-criterion}, every perfect-Delaunay facet inequality is valid for $\mathcal C_n$.  Proposition~\ref{prop:perfect-facets-describe-hull} therefore gives $\mathcal C_n\subseteq\MIQ_n$.  The reverse inclusion holds by \eqref{eq:universal-chain}, so $\mathcal C_n=\MIQ_n$. Conversely, suppose that $\mathcal C_n=\MIQ_n$.  Every facet-defining inequality of $\MIQ_n$ is then valid for $\mathcal C_n$.  Corollary~\ref{cor:depth-criterion} implies that every perfect Delaunay polyhedron in $\Z^n$ is shallow.
\end{proof}

In the next section, we will use the classifications of perfect Delaunay polyhedra and their depths to determine the dimensions in which $\mathcal C_n=\MIQ_n$.

%% file: sections/05_exact_hull.tex
\section{Exact recovery of \texorpdfstring{$\MIQ_n$}{MIQ\_n}}
\label{sec:low-dimensional}

Our results in the previous sections reduce the comparison of $\mathrm{MIQ}_n$, $\mathcal C_n$, and $\mathcal S_n$ to the geometry of perfect Delaunay polyhedra. By Theorem~\ref{thm:facet-correspondence}, these polyhedra correspond to the facets of $\mathrm{MIQ}_n$, while Corollary~\ref{cor:depth-criterion} determines whether the associated facet-defining inequality is obtained by one round of CG cuts. We now use the known classification of perfect Delaunay polyhedra in low dimensions to make these results concrete. In particular, we identify the perfect Delaunay polyhedra that occur up to dimension seven, compute their depths, and then use these calculations to determine exactly when $\mathcal C_n$ and $\mathcal S_n$ coincide with $\MIQ_n$. Specifically, we prove that $\MIQ_n=\mathcal C_n$ for $1\le n \le 6$, $\MIQ_n=\mathcal S_n$ for $1\le n\le 5$, and $\mathcal C_n\subsetneq\mathcal S_n$ for every $n\ge6$. We further prove the strict inclusions in higher dimensions using the results established in Section~\ref{sec:depth}.

Throughout this section, we use \emph{type} to refer to an equivalence class
of perfect Delaunay polyhedra in a fixed dimension under integer-affine
bijections.  Perfect Delaunay polytopes have been classified through
dimension seven: Erdahl treats dimensions at most five
\cite[Theorem~5.1 and Corollary~5.2]{Erdahl1992}, and Dutour Sikiri\'c treats
dimensions six \cite[Theorem~2]{Dutour2004} and seven
\cite[Theorem~1.1]{DutourSikiric2017}. For the seven-dimensional polytopes,
we use the quadratic-form matrices given by Dutour Sikiri\'c, Erdahl, and
Rybnikov \cite[Sections~7.1 and~7.2]{DutourErdahlRybnikov2007}, with
the defining quadratics specified in Appendix~\ref{app:low-dimensional}.
As in Section~\ref{sec:delaunay}, we use the names of these polytopes for
their lattice-point sets when discussing integer-affine equivalence,
computing depth, or taking products with $\Z^k$.

\begin{theorem}[Depth classification up to dimension seven]
\label{thm:depth-classification}
Up to integer-affine equivalence, the perfect Delaunay polyhedra in $\Z^n$ for $1\le n\le7$, together with their depths and classification as shallow or deep, are those in Table~\ref{tab:depths}. In particular, $ER_7$ is the only deep type up to dimension seven, with depth $43/32$.
\end{theorem}

\begin{table}[htbp]
\centering
\caption{Perfect-Delaunay types and depths up to dimension seven.}
\label{tab:depths}
\begin{tabular}{clcc}
\toprule
Dimension & Type & Depth & Status\\
\midrule
$1\le n\le7$ & $\{0,1\}\times\Z^{n-1}$ (primitive split) & $1/8$ & shallow\\
$6$ & $2_{21}$ & $2/3$ & shallow\\
$7$ & $2_{21}\times\Z$ & $2/3$ & shallow\\
$7$ & $3_{21}$ & $3/4$ & shallow\\
$7$ & $ER_7$ & $43/32$ & deep\\
\bottomrule
\end{tabular}
\end{table}

\begin{proof}
By Theorem~\ref{thm:perfect-delaunay-structure},
every perfect Delaunay polyhedron $D\subseteq\Z^n$ has, after an
integer-affine change of coordinates, the form
\[
 D=D_0\times\Z^{n-k},\qquad 1\le k\le n,
\]
where $D_0\subseteq\Z^k$ is the lattice-point set of a $k$-dimensional
perfect Delaunay polytope.  Proposition~\ref{prop:depth-properties} gives
$\delta(D)=\delta(D_0)$.  Thus it suffices to classify these polytopes
and compute their depths.

The cited classifications of perfect Delaunay polytopes provide the following types in dimensions up to seven: the type
$\{0,1\}$ in dimension one, no types in dimensions two to five, the Gosset type $2_{21}$
in dimension six, and $3_{21}$ and $ER_7$ in dimension seven.
The products $\{0,1\}\times\Z^{n-1}$ represent the primitive-split
Delaunay polyhedra $D_{v,s}$ from \eqref{eq:split-Delaunay}.
Up to dimension seven, the only further product with a positive-dimensional
integer lattice is $2_{21}\times\Z$.  These types are summarized in Table~\ref{tab:depths}.

For each integer-valued defining quadratic $q\in\mathcal Q_n$, the depth
is $-\min_{y\in\R^n}q(y)/\Delta(q)$ by Definition~\ref{def:depth} and
\eqref{eq:primitive-content}.  Appendix~\ref{app:low-dimensional}
presents the coefficient and positive-definiteness checks, so we give
only the resulting minima and primitive normalizations here.

For the type $\{0,1\}$, take $v=1$ and $s=0$ in
\eqref{eq:split-polynomial}.  Completing the square gives
\[
 \widehat q_{1,0}(z)=\binom{z}{2}
 =\frac12\left(z-\frac12\right)^2-\frac18.
\]
Its nonconstant Newton coefficient vector is $(0,1)$, so
$\Delta(\widehat q_{1,0})=1$ and $\delta(\{0,1\})=1/8$.

Burer and Letchford \cite[Section~5.3, (10)]{BurerLetchford2014} exhibited a
bounded facet of $\MIQ_6$ containing $27$ extreme points at equality. They
observed, without giving a formal proof, that the convex hull of the
corresponding $x$-vectors is an affine image of the six-dimensional Gosset
polytope $2_{21}$. For this type, we specialize their lifted facet inequality
to obtain the defining quadratic
\begin{equation}\label{eq:G6-generator}
q_{G_6}(z) =\sum_{i=1}^{6}(z_i^2-z_i)+2z_6 +z_1z_2+z_2z_3+z_3z_4+z_4z_5+z_5z_1 -z_6\sum_{i=1}^{5}z_i.
\end{equation}
Because the cited inequality is facet-defining for $\MIQ_6$,
Theorem~\ref{thm:facet-correspondence} shows that $Z(q_{G_6})$ is perfect.
Appendix~\ref{app:low-dimensional} verifies that its quadratic-form matrix
is positive definite, so $Z(q_{G_6})$ is finite.
The six-dimensional classification \cite[Theorem~2]{Dutour2004} therefore
identifies $Z(q_{G_6})$ as an integer-affine image of $2_{21}$. In addition,
Appendix~\ref{app:low-dimensional} gives
\[
\min_{\R^6}q_{G_6}=-\frac{2}{3}, \qquad \Delta(q_{G_6})=1, \qquad \widehat q_{G_6}=q_{G_6}, \qquad \delta(2_{21})=\frac{2}{3}.
\]

For $3_{21}$ and $ER_7$, use the defining quadratics $q_{G_7}$ and $q_{ER}$ in Appendix~\ref{app:low-dimensional}, where we give the explicit matrices and verify that both nonconstant Newton coefficient vectors have greatest common divisor two. Hence,
\[
\begin{aligned}
\min_{\R^7}q_{G_7}&=-\frac{3}{2}, &\widehat q_{G_7}&=\frac{1}{2}q_{G_7}, &\delta(3_{21})&=\frac{3}{4},\\
\min_{\R^7}q_{ER}&=-\frac{43}{16}, &\widehat q_{ER}&=\frac{1}{2}q_{ER}, &\delta(ER_7)&=\frac{43}{32}>1.
\end{aligned}
\]

These calculations and the depth equality above give all values and
classifications as shallow or deep in Table~\ref{tab:depths}.
\end{proof}

The depth calculations now determine the facet-defining inequalities obtained from $\mathcal P_n$ by CG rounding. By Proposition~\ref{prop:facet-normal-rounding}, rounding the inequalities $L_{\widehat q_{G_6}}\ge-2/3$ over $\mathcal P_6$ and $L_{\widehat q_{G_7}}\ge-3/4$ over $\mathcal P_7$, using the nonconstant Newton coefficient vectors of these primitive integer-valued generators, gives the corresponding facet-defining inequalities with right-hand side zero. For $ER_7$, rounding gives only $L_{\widehat q_{ER}}\ge-1$. Corollary~\ref{cor:depth-criterion} further shows that $L_{\widehat q_{ER}}\ge0$ is not valid for $\mathcal C_7$. We next combine these facts with the description of $\MIQ_n$ in Proposition~\ref{prop:perfect-facets-describe-hull} to determine exactly when the two closures describe the integer quadratic moment hull.

The inclusions in \eqref{eq:universal-chain} hold in every dimension. We first establish the equality cases. The classification up to dimension six, together with the description of $\MIQ_n$ in Proposition~\ref{prop:perfect-facets-describe-hull}, is sufficient for this step.

\begin{lemma}
\label{lem:low-dimensional-equalities}
For $1\le n\le5$, $\MIQ_n=\mathcal C_n=\mathcal S_n$. Moreover, $\MIQ_6=\mathcal C_6$. 
\end{lemma}

\begin{proof}
For $1\le n\le5$, Theorem~\ref{thm:depth-classification} shows that every perfect Delaunay polyhedron is integer-affinely equivalent to $\{0,1\}\times\Z^{n-1}$. These are exactly the primitive-split Delaunay polyhedra $D_{v,s}$. Theorem~\ref{thm:facet-correspondence} therefore identifies every facet-defining inequality of $\MIQ_n$ as a primitive split inequality, up to positive scaling. Because these inequalities describe $\MIQ_n$ by Proposition~\ref{prop:perfect-facets-describe-hull}, the inclusions in \eqref{eq:universal-chain} give $\MIQ_n=\mathcal C_n=\mathcal S_n$.

In dimension six, Theorem~\ref{thm:depth-classification} shows that every perfect Delaunay polyhedron is integer-affinely equivalent to $\{0,1\}\times\Z^5$ or to $2_{21}$. Both are shallow. Corollary~\ref{cor:depth-exact-recovery} therefore gives $\MIQ_6=\mathcal C_6$.
\end{proof}

To prove $\mathcal C_n\subsetneq\mathcal S_n$ in higher dimensions, we next exhibit a facet-defining inequality of $\MIQ_n$ that is valid for $\mathcal C_n$ but not for $\mathcal S_n$. The inequalities associated with $2_{21}$ and its products with integer lattices give such examples in every dimension $n\ge6$.

\begin{lemma}
\label{lem:CG-strictly-stronger-than-splits}
For every $n\ge6$, $\mathcal C_n\subsetneq\mathcal S_n$.
\end{lemma}

\begin{proof}
Fix $n\ge6$ and define
\[
\widehat q_{G_6}^{(n)}(z) :=\widehat q_{G_6}(z_1,\ldots,z_6),\qquad z\in\R^n.
\]
For $n=6$, this is the original primitive integer-valued generator. For $n>6$, its lattice-zero set is $Z(q_{G_6})\times\Z^{n-6}$, which is perfect and has depth $2/3$ by Proposition~\ref{prop:depth-properties}. The nonconstant Newton coefficient vector remains primitive because introducing additional variables only inserts zero coefficients. Thus, in every dimension $n\ge6$, the inequality $L_{\widehat q_{G_6}^{(n)}}\ge0$ is facet-defining for $\MIQ_n$ and valid for $\mathcal C_n$ by Theorem~\ref{thm:facet-correspondence} and Corollary~\ref{cor:depth-criterion}.

To see that this inequality is not valid for $\mathcal S_n$, apply Theorem~\ref{thm:rigidity-depth} to the family
\[
\mathcal F:=\{\widehat q_{v,s}:v\in\Z^n\text{ primitive},\ s\in\Z\}.
\]
Every polynomial in $\mathcal F$ is integer-valued and has continuous minimum $-1/8>-1$. Moreover, $\mathcal K_{\mathcal F}=\mathcal S_n$ because $\mathcal S_n\subseteq\mathcal P_n$ by Proposition~\ref{prop:splits-imply-psd}. The quadratic-form matrix of $\widehat q_{G_6}^{(n)}$ has rank six, whereas that of $\widehat q_{v,s}$ is $vv^{\top}/2$ and has rank one. Hence, $\mathcal F$ contains no positive integral multiple of $\widehat q_{G_6}^{(n)}$. Theorem~\ref{thm:rigidity-depth} therefore shows that its facet-defining inequality is not valid for $\mathcal S_n$. Because the same inequality is valid for $\mathcal C_n$, we obtain $\mathcal C_n\subsetneq\mathcal S_n$ for every $n\ge6$.
\end{proof}

We now show that $\mathcal C_n$ itself fails to equal $\MIQ_n$ for $n\ge7$, using the facet-defining inequalities associated with $ER_7$ and its products with integer lattices.

\begin{lemma}
\label{lem:CG-not-exact-from-seven}
For every $n\ge7$, $\MIQ_n\subsetneq\mathcal C_n$.
\end{lemma}

\begin{proof}
Fix $n\ge7$ and define
\[
\widehat q_{ER}^{(n)}(z) :=\widehat q_{ER}(z_1,\ldots,z_7),\qquad z\in\R^n.
\]
For $n=7$, its lattice-zero set is the given representative of $ER_7$ and has depth $43/32$. For $n>7$, its lattice-zero set is $Z(q_{ER})\times\Z^{n-7}$, which remains perfect and has the same depth by Proposition~\ref{prop:depth-properties}. Thus, $L_{\widehat q_{ER}^{(n)}}\ge0$ is a facet-defining inequality of $\MIQ_n$ that is not valid for $\mathcal C_n$, by Theorem~\ref{thm:facet-correspondence} and Corollary~\ref{cor:depth-criterion}. This proves the strict inclusion for every $n\ge7$.
\end{proof}

With these three lemmas, we give the complete comparison of $\MIQ_n$ with the two closures $\mathcal C_n$ and $\mathcal S_n$. We also identify the facet-defining inequalities that remain valid for $\mathcal C_7$ with Table~\ref{tab:depths}. In turn, we determine the dimensions in which $\MIQ_n$ is described by the split inequalities or by the rank-one CG closure of $\mathcal P_n$. 

\begin{theorem}
\label{thm:closure-comparison}
For every $n\ge1$,
\[
\begin{array}{rll}
1\le n\le5: & \MIQ_n=\mathcal C_n=\mathcal S_n, &\\[1mm]
n=6:         & \MIQ_6=\mathcal C_6\subsetneq\mathcal S_6, &\\[1mm]
n\ge7:      & \MIQ_n\subsetneq\mathcal C_n\subsetneq\mathcal S_n. &
\end{array}
\]
Moreover, $\mathcal S_n\subsetneq\mathcal P_n$ for $n \ge 1$. For every perfect Delaunay polyhedron $D\subseteq\Z^7$, $L_{\widehat q_D}\ge0$ is valid for $\mathcal C_7$ if and only if $D$ is not integer-affinely equivalent to $ER_7$. 
\end{theorem}

\begin{proof}
Lemma~\ref{lem:low-dimensional-equalities} gives the equality cases. Lemmas~\ref{lem:CG-strictly-stronger-than-splits} and~\ref{lem:CG-not-exact-from-seven} give the strict inclusions between $\MIQ_n$, $\mathcal C_n$, and $\mathcal S_n$. Proposition~\ref{prop:splits-imply-psd} gives $\mathcal S_n\subsetneq\mathcal P_n$. Finally, Corollary~\ref{cor:depth-criterion} says that $L_{\widehat q_D}\ge0$ is valid for $\mathcal C_7$ exactly when $\delta(D)<1$. By Theorem~\ref{thm:depth-classification}, this holds exactly when $D$ is not integer-affinely equivalent to $ER_7$.
\end{proof}

In dimension seven, adding the facet inequalities associated with all integer-affine images of $ER_7$ to $\mathcal C_7$ gives the full moment hull. Let $T(z)=Uz+t$ range over all integer-affine bijections of $\R^7$, where $U\in\GL_7(\Z)$ and $t\in\Z^7$. Then
\[
\MIQ_7 =\mathcal C_7\cap \bigcap_T \{W\in\mathcal H_7:L_{q_{ER}\circ T^{-1}}(W)\ge0\}.
\]
The added inequalities give all facet inequalities of $\MIQ_7$ associated with integer-affine images of $ER_7$. Every other perfect-Delaunay facet inequality is valid for $\mathcal C_7$ by Theorem~\ref{thm:closure-comparison}, so Proposition~\ref{prop:perfect-facets-describe-hull} gives the equality.

%% file: sections/06_complexity.tex
\section{Complexity of membership and separation}
\label{sec:complexity}

We now study the complexity of membership and exact
separation for the rank-one CG closure $\mathcal C_n$, when the
dimension $n$ is part of the input. Given a rational point
$\bar\eta\in\mathbb Q^{d_n}$, we consider the decision problems
of determining whether $\bar\eta\notin\mathcal C_n$
($\CGViol$) or $\bar\eta\in\mathcal C_n$ ($\CGMem$).
Exact separation means either certifying membership or finding
a violated rank-one CG inequality. We first show that every
rational point outside $\mathcal C_n$ admits a polynomial-size
certificate of nonmembership that can be verified in polynomial
time. We then establish strong NP-completeness of $\CGViol$
and strong NP-hardness of exact separation over both
$\mathcal C_n$ and the split-inequality closure $\mathcal S_n$.

\subsection{Polynomial-size certificates}
\label{subsec:polynomial-certificates}
\label{subsec:zero-extension}

Throughout this section, we measure input size by the total binary encoding
length of $(n,\bar\eta)$, with rational coordinates represented as reduced
fractions.

By Lemma~\ref{lem:CG-polynomial-dictionary}, $\bar\eta\notin\mathcal C_n$ if and only if
there exists $p\in\operatorname{Int}_2(\Z^n)$
such that
\begin{equation}\label{eq:certificate-equivalence}
-1<\min_{y\in\mathbb R^n}p(y)\le 0,
\quad L_p(\bar\eta)<0.
\end{equation}
We call such a polynomial a \emph{cut-polynomial
certificate} of nonmembership. We encode $p$ by its
ordered vector of Newton coefficients, as defined in
Section~\ref{sec:coordinates}.

Given a proposed certificate
$p(z)=\alpha_p+c_p^\top\eta(z)$, we can verify
its validity in polynomial time. First, we check that
its Newton coefficients are integers. Using~\eqref{eq:coefficient-dictionary},
write
\[
p(z)=z^\top A_pz+2b_p^\top z+\alpha_p.
\]
The continuous minimum of $p$ is finite if and only if
\[
A_p\succeq 0
\quad\text{and}\quad
b_p\in\operatorname{range}(A_p).
\]
Both conditions can be checked by exact rational linear
algebra. When they hold, we solve $A_py=-b_p$ and obtain
\[
\min_{z\in\mathbb R^n}p(z)
=\alpha_p+b_p^\top y,
\]
which is independent of the chosen solution. We then
check whether this minimum belongs to $(-1,0]$ and
evaluate
\[
L_p(\bar\eta)=\alpha_p+c_p^\top\bar\eta
\]
to determine whether the inequality is violated.
All these operations have polynomial bit complexity
in the combined encoding lengths of the input and
the certificate.

The construction below uses zero extension, by which we mean enlarging a moment matrix by adjoining zero rows and columns. The next lemma shows that this operation preserves membership in the rank-one CG
closure.

Fix integers $1\le m<n$. Using the coordinate bijection in Lemma~\ref{lem:integral-moment-coordinates}, we define the linear map $\iota:\R^{d_m}\to\R^{d_n}$ by retaining all integral moment coordinates involving only the first $m$ variables and setting every coordinate involving one of the last $n-m$ variables to zero. Equivalently, after ordering the indices so that the new variables come last,
\[
W(\iota(\eta))=
\begin{bmatrix}
W(\eta)&0\\
0&0_{n-m}
\end{bmatrix}
\qquad(\eta\in\R^{d_m}),
\]
where $0_{n-m}$ is the $(n-m)\times(n-m)$ zero matrix and the off-diagonal blocks are zero matrices of the required sizes. The block formula uniquely determines $\iota(\eta)$ because $\eta\mapsto W(\eta)$ is a bijection. The following lemma establishes the CG-closure membership equivalence.

\begin{lemma}
\label{lem:zero-extension}
For every $\eta\in\R^{d_m}$,
\[
\eta\in\mathcal C_m\quad\Longleftrightarrow\quad\iota(\eta)\in\mathcal C_n.
\]
\end{lemma}

\begin{proof}
We first record an evaluation identity for lifting certificates. For $q\in\mathcal Q_m$, define $\widetilde q\in\mathcal Q_n$ by
\[
\widetilde q(z):=q(z_1,\ldots,z_m),\qquad z\in\R^n.
\]
With the variables ordered as above, the coefficient matrix of $\widetilde q$ is
\[
M_{\widetilde q}=
\begin{bmatrix}
M_q&0\\
0&0_{n-m}
\end{bmatrix}.
\]
For every $\eta\in\R^{d_m}$, the definition of $\iota$ gives
\[
L_{\widetilde q}(\iota(\eta)) =\ip{M_{\widetilde q}}{W(\iota(\eta))} =\ip{M_q}{W(\eta)} =L_q(\eta).
\]
For the membership equivalence, first suppose that $\eta\in\mathcal C_m$. Consider any $c\in\Z^{d_n}$ satisfying $\beta_{\mathcal P_n}(c)>-\infty$. From this vector, form $c'\in\Z^{d_m}$ by retaining the coefficients of the nonconstant Newton basis functions involving only the first $m$ variables. Thus,
\[
{c'}^{\top}\eta(y)=c^{\top}\eta((y,0))\qquad(y\in\R^m),
\]
where $(y,0)\in\R^n$ has its last $n-m$ coordinates equal to zero. By \eqref{eq:beta-formula},
\[
\begin{aligned}
\beta_{\mathcal P_m}(c') &=\inf_{y\in\R^m}c^{\top}\eta((y,0))\\
&\ge\inf_{z\in\R^n}c^{\top}\eta(z) =\beta_{\mathcal P_n}(c)>-\infty.
\end{aligned}
\]
Therefore,
\[
c^{\top}\iota(\eta) ={c'}^{\top}\eta \ge\left\lceil\beta_{\mathcal P_m}(c')\right\rceil \ge\left\lceil\beta_{\mathcal P_n}(c)\right\rceil,
\]
where the first inequality is the defining CG inequality for $c'$; when $c'=0$, it is the trivial inequality $0\ge0$. Because $c$ was arbitrary, $\iota(\eta)\in\mathcal C_n$. On the other hand, suppose that $\iota(\eta)\in\mathcal C_n$. Consider any $c\in\Z^{d_m}$ satisfying $\beta_{\mathcal P_m}(c)>-\infty$. Form $\widetilde c\in\Z^{d_n}$ by retaining these Newton coefficients and assigning coefficient zero to every basis function involving an additional variable. We thereby obtain a quadratic that is independent of the additional variables, so \eqref{eq:beta-formula} gives
\[
\beta_{\mathcal P_n}(\widetilde c)=\beta_{\mathcal P_m}(c), \qquad \widetilde c^{\top}\iota(\eta)=c^{\top}\eta.
\]
Apply the defining CG inequality for $\widetilde c$ at $\iota(\eta)$ and use these identities to obtain $c^{\top}\eta\ge\lceil\beta_{\mathcal P_m}(c)\rceil$. Because $c$ was arbitrary, $\eta$ satisfies every defining inequality of $\mathcal C_m$.

\end{proof}

Now we show that a certificate of polynomial
encoding length always exists.

\begin{theorem}[Polynomial-size cut-polynomial certificates]
\label{thm:polynomial-certificates}
Every rational point $\bar\eta\notin\mathcal C_n$ admits a cut-polynomial
certificate whose binary encoding length is bounded by a polynomial in the
input length of $(n,\bar\eta)$.  Consequently,
$\CGViol\in\mathrm{NP}$.
\end{theorem}

\begin{proof}
Let
\[
\bar W=
\begin{bmatrix}
1&\bar x^\top\\
\bar x&\bar X
\end{bmatrix},
\qquad
\Sigma:=\bar X-\bar x\bar x^\top.
\]
By~\eqref{eq:certificate-equivalence}, every
$\bar\eta\notin\mathcal C_n$ admits a violated
cut-polynomial certificate. We show that such a
certificate can be chosen with polynomial binary
encoding length by considering three cases.

\medskip
\noindent\emph{Indefinite moment matrix.}
Suppose first that $\bar W$ is indefinite.  Since $\bar W$ has a negative eigenvalue,
scaling and rounding a corresponding eigenvector
yields an integer vector
$v=(v_0,v')\in\Z^{n+1}$ of polynomial encoding
length such that
\[
v^\top\bar Wv<0.
\]
Appendix~\ref{subsec:encoding-indefinite} shows the scaling factor and provides the bit-length bound.   Necessarily $v'\ne0$, because
$\bar W_{00}=1$.  The affine square
\[
 p_v(z):=\left(v_0+{v'}^{\top}z\right)^2
\]
is integer-valued, has continuous minimum zero, and satisfies
\[
 L_{p_v}(\bar\eta)=v^{\top}\bar Wv<0.
\]
Its Newton coefficients are sums and pairwise products of entries of $v$, so
it is a polynomial-size certificate.  

\medskip
\noindent\emph{Positive-definite covariance.}
Suppose $\bar W\succeq0$ and $\Sigma\succ0$. Choose any certificate $p$ from
\eqref{eq:certificate-equivalence}, and write
\[
 p(z)=z^{\top}Az+2b^{\top}z+\alpha,
 \qquad A\succeq0,
 \qquad b\in\range(A).
\]
Let $\min p=-\delta$, where
$0\le\delta<1$.  A violated certificate cannot have $\delta=0$,
because then
\[
 L_p(\bar\eta)=p(\bar x)+\ip{A}{\Sigma}\ge0.
\]
Thus $0<\delta<1$.  Choose $y$ satisfying $Ay=-b$. Completing the square gives
\begin{equation}\label{eq:positive-covariance-master}
 L_p(\bar\eta)
 =-\delta+(\bar x-y)^{\top}A(\bar x-y)+\ip{A}{\Sigma}.
\end{equation}
Since both terms following $-\delta$ are
nonnegative, the violation implies
\begin{equation}\label{eq:positive-covariance-trace}
 \ip{A}{\Sigma}<1.
\end{equation}
Let $\lambda=\lambda_{\min}(\Sigma)>0$, then
$\lambda\operatorname{tr}(A)\le\ip{A}{\Sigma}$, and hence
\[
 \|A\|_2\le\operatorname{tr}(A)<\lambda^{-1}.
\]
Because $\Sigma$ is rational and positive definite, the binary length of a
valid upper bound on $\lambda^{-1}$ is polynomial in the input length.
The identity~\eqref{eq:positive-covariance-master} then bounds $b$ and $\alpha$,
and hence all Newton coefficients of $p$, by numbers of polynomial binary
length.  The quantitative estimates are given in
Appendix~\ref{subsec:encoding-positive-definite}.

\medskip
\noindent\emph{Singular covariance.} It remains to consider $\bar W\succeq0$ with
singular covariance $\Sigma$. Let
$r=\rank(\Sigma)<n$.  In this case, the bound
$\langle A,\Sigma\rangle<1$ still holds, but no longer controls $A$ along directions in
$\ker(\Sigma)$. We therefore find an integer basis of $\ker(\Sigma)\cap\Z^n$ and handle these directions separately.

Choose a positive integer
$D_\Sigma$ that clears the denominators of
$\Sigma$, and set
\[
S=D_\Sigma\Sigma\in\Z^{n\times n}.
\]

When $r>0$, since $\ker(S)=\ker(\Sigma)$, we can apply the full-row-rank extension of the
Kannan--Bachem Hermite-normal-form algorithm
\cite[Section~2, especially Theorem~2 and p.~504]{KannanBachem} to $r$
linearly independent rows of $S$ to obtain a
unimodular matrix\[
Q=[Q_1\ Z]\in\mathrm{GL}_n(\Z),
\]
where $Q_1\in\Z^{n\times r}$ consists of the first
$r$ columns of $Q$, and the columns of
$Z\in\Z^{n\times(n-r)}$ form a lattice basis of
$\ker(\Sigma)\cap\Z^n$. When $r=0$, we take
$Q=I_n$. The matrix $Q$ has polynomial encoding
length; Appendix~\ref{subsec:encoding-singular}
provides the verification and the bit bounds.

Suppose first that $Z^\top\bar x\notin\Z^{n-r}$.
Then some column $v$ of $Z$ satisfies
$v^\top\bar x\notin\Z$. Since $Q$ is unimodular,
$v$ is primitive. Let
$s=\lfloor v^\top\bar x\rfloor$ and consider
the split quadratic
\[
\widehat q_{v,s}(z)
=\frac12(v^\top z-s)(v^\top z-s-1).
\]
Because $\Sigma v=0$, we have
\[
L_{\widehat q_{v,s}}(\bar\eta)
=\frac12\theta(\theta-1)<0,
\qquad
\theta:=v^\top\bar x-s\in(0,1).
\]
Thus, $\widehat q_{v,s}$ is a violated
cut-polynomial certificate. Its coefficients
have polynomial encoding length because
$v$ and $s$ do.

Otherwise, $Z^\top\bar x\in\Z^{n-r}$.
In this case, we can eliminate the zero-variance
directions by an integer-affine transformation.
Let
\[
U=Q^{-\top},
\qquad
t=U
\begin{bmatrix}
0\\ Z^\top\bar x
\end{bmatrix}
\in\Z^n.
\]
Under the change of variables $z=t+Uy$,
the
transformed mean is
\[
\bar y
=U^{-1}(\bar x-t)
=Q^\top\bar x-
\begin{bmatrix}0\\Z^\top\bar x\end{bmatrix}
=
\begin{bmatrix}Q_1^\top\bar x\\0\end{bmatrix}.
\]
Similarly, since $\Sigma Z=0$, the transformed
covariance is
\[
\Sigma'
=U^{-1}\Sigma U^{-\top}
=Q^\top\Sigma Q
=
\begin{bmatrix}
Q_1^\top\Sigma Q_1&0\\
0&0
\end{bmatrix}.
\]
Defining $\bar y_r:=Q_1^\top\bar x, \Sigma_r:=Q_1^\top\Sigma Q_1$,
we obtain
\[
\bar y=
\begin{bmatrix}\bar y_r\\0\end{bmatrix},
\qquad
\Sigma'=
\begin{bmatrix}\Sigma_r&0\\0&0\end{bmatrix},
\quad \Sigma_r\succ0.
\]

If $r=0$, the transformed point is $W_0\in
\mathcal C_n$, contradicting
$\bar\eta\notin\mathcal C_n$. Hence $r\ge1$. Since the columns of $Z$ span $\ker(\Sigma)$
and $Q$ is invertible, we have $\Sigma_r\succ0$.
Thus, the transformed moment point is the exact
zero extension of an $r$-dimensional point
$\bar\eta_r$ with positive-definite covariance. 

By integer-affine invariance and
Lemma~\ref{lem:zero-extension},
$\bar\eta_r\notin\mathcal C_r$. The positive-definite case therefore gives a polynomial-size cut-polynomial certificate
$p_r$ for $\bar\eta_r$. Extending $p_r$ by zero variables and applying
the inverse coordinate transformation yields
a certificate for $\bar\eta$. Since the
transformation and its inverse have polynomial
encoding length, the resulting Newton
coefficients also have polynomial encoding
length. The details are given in
Appendix~\ref{subsec:encoding-singular}.

In every case, $\bar\eta\notin\mathcal C_n$
admits a polynomial-size cut-polynomial
certificate.
\end{proof}

\subsection{Strong NP-completeness and consequences}
\label{subsec:strong-hardness}

Having shown that $\CGViol$ belongs to NP, we now establish its strong
NP-hardness by a reduction from a gap version of the nearest-codeword
problem. The reduction produces rational input points with numerator and
denominator magnitudes polynomially bounded in $n$.

The nearest-codeword problem seeks a vector $u$ that minimizes the number of violated equations in a linear system over a finite field. We use a gap version that distinguishes systems in which some assignment violates few equations from systems in which every assignment violates many equations. The following hardness result is a consequence of H{\aa}stad's inapproximability theorem~\cite{Hastad01}. The notations $\F_{11}$ and $\Z_+$  denote the field of integers modulo 11 and nonnegative integers, respectively. 

\begin{theorem}[Fixed-field gap nearest codeword; consequence of
H{\aa}stad~\cite{Hastad01}]
\label{thm:gap-nearest-codeword}
Given
\[
 A\in\F_{11}^{N\times r},
 \qquad
 t\in\F_{11}^{N},
 \qquad
 k\in\Z_+,
\]
define
\[
 D_H:=\min_{u\in\F_{11}^r}\|Au-t\|_0.
\]
It is NP-hard to distinguish instances satisfying
\[
 D_H\le k
 \qquad\text{from those satisfying}\qquad
 D_H>100k.
\]
\end{theorem}

\begin{proof}
The proof of H{\aa}stad's theorem on \textsc{Max-E3-Lin} over finite Abelian
groups---specifically, its completeness and soundness lemmas---gives the
following statement for the additive group of $\F_{11}$
\cite[Theorem~5.9 and Lemmas~5.10--5.11]{Hastad01}.  For every fixed
\[
 \rho\in\left(0,\frac{1-1/11}{2}\right),
\]
it is NP-hard to distinguish a system of $N_0$ linear equations over
$\F_{11}$ for which some assignment satisfies at least $(1-\rho)N_0$
equations from one for which every assignment satisfies at most
$(1/11+\rho)N_0$ equations.

Set $\rho:=1/220$ and repeat every equation exactly $220$ times.  Let
$A\in\F_{11}^{N\times r}$ be the coefficient matrix of the repeated
system, let $t\in\F_{11}^N$ be its right-hand side, and set
\[
 N:=220N_0,
 \qquad
 k:=N_0.
\]
For every assignment $u$, the Hamming weight $\|Au-t\|_0$ is exactly the
number of repeated equations violated by $u$.

In the YES case, some assignment violates at most $N_0/220$ of the original
equations.  After repetition, it therefore violates at most
\[
 220\frac{N_0}{220}=N_0=k
\]
equations, and hence $D_H\le k$.  In the NO case, every assignment violates
at least
\[
 \left(1-\frac1{11}-\frac1{220}\right)N_0
 =\frac{199}{220}N_0
\]
of the original equations.  Consequently,
\[
 D_H\ge220\frac{199}{220}N_0
 =199N_0
 =199k
 >100k.
\]

If $A$ is not of full column rank, compute a basis of its column space over
$\F_{11}$ and replace the columns of $A$ by that basis.  This preserves the
linear code $\{Au:u\in\F_{11}^r\}$, and therefore preserves $D_H$, while
making $A$ full column rank.  The repetition and column-basis operations are
deterministic and polynomial-time, and the construction satisfies
$N=220k$.
\end{proof}

The reduction below uses the instances constructed above, with $A$ of full
column rank and $N=220k$. We next construct an integer least-squares
problem whose minimum is a fixed multiple of $D_H$. After minimizing over auxiliary integer variables, every violated equation contributes the same amount to the squared error.
 Choose integer representatives
\[
 \widehat A\in\{0,\ldots,10\}^{N\times r},
 \qquad
 \widehat t\in\{0,\ldots,10\}^N.
\]
Let $M:=10N$, index its rows by pairs
$(j,a)$ with $j=1,\ldots,N$ and $a=1,\ldots,10$, and put
$m:=r+M$.  For $h=(u,\zeta)\in\Z^r\times\Z^M$, define
$B\in\Z^{M\times m}$ and $\tau\in\Z^M$ by
\[
 (Bh)_{j,a}:=a(\widehat Au)_j+11\zeta_{j,a},
 \qquad
 \tau_{j,a}:=a\widehat t_j.
\]

\begin{lemma}[Scalar invariance]
\label{lem:scalar-invariance}
For every integer $d$ with $11\nmid d$,
\[
 \min_{h\in\Z^m}\|Bh-d\tau\|_2^2=110D_H.
\]
\end{lemma}

\begin{proof}
Fix $d$ and first fix the integer vector $u$.  For coordinate $j$, let
$r_j:=(Au-dt)_j\in\F_{11}$.  Minimizing over the independent variables
$\zeta_{j,a}$ gives
\[
 \sum_{a=1}^{10}\operatorname{dist}(ar_j,11\Z)^2.
\]
The contribution is zero when $r_j=0$.  When $r_j\ne0$, multiplication by
$r_j$ permutes the nonzero residues modulo $11$, whose centered
representatives are
\[
 \pm1,\ \pm2,\ \pm3,\ \pm4,\ \pm5.
\]
The contribution is therefore
\[
 2(1^2+2^2+3^2+4^2+5^2)=110.
\]
After minimizing over $\zeta$,
\[
 \min_{\zeta}\|B(u,\zeta)-d\tau\|_2^2
 =110\|Au-dt\|_0.
\]
Changing $u$ by an element of $11\Z^r$ can be absorbed into $\zeta$, so
minimization over integer $u$ is minimization over $\F_{11}^r$.
Multiplication by $d\ne0$ bijects $\F_{11}^r$ and preserves Hamming weight,
so
\[
 \min_{u\in\F_{11}^r}\|Au-dt\|_0
 =\min_{u\in\F_{11}^r}\|Au-t\|_0=D_H.
\]
\end{proof}

Define the positive-semidefinite Gram matrix
\[
 G:=\begin{bmatrix}-\tau&B\end{bmatrix}^{\top}
    \begin{bmatrix}-\tau&B\end{bmatrix}\succeq0.
\]
Thus, for $d\in\Z$ and $h\in\Z^m$,
\begin{equation}\label{eq:hard-Gram-distance}
 \begin{bmatrix}d\\h\end{bmatrix}^{\top}
 G
 \begin{bmatrix}d\\h\end{bmatrix}
 =\|Bh-d\tau\|_2^2.
\end{equation}
Set
\[
 \varepsilon:=\frac1{44},
 \qquad
 T:=220k,
 \qquad
 n:=m+1.
\]
The first decision coordinate is the anchor, with unit vector $e_1\in\R^n$.
Define
\begin{equation}\label{eq:hard-point}
 \bar x:=\frac12e_1,
 \qquad
 \Sigma:=\left(\frac14-\varepsilon\right)e_1e_1^{\top}
          +\frac{\varepsilon}{T}G,
 \qquad
 \bar X:=\bar x\bar x^{\top}+\Sigma.
\end{equation}
Because $G\succeq0$ and $1/4-\varepsilon=5/22>0$, the Schur complement
shows that the associated rational point $\bar\eta$ belongs to
$\mathcal P_n$.

\begin{lemma}[Extraction from an arbitrary violated cut]
\label{lem:arbitrary-cut-extraction}
If $\bar\eta\notin\mathcal C_n$, then there exist
$d\in\{1,\ldots,7\}$ and $h\in\Z^m$ such that
\[
 \begin{bmatrix}d\\h\end{bmatrix}^{\top}
 G
 \begin{bmatrix}d\\h\end{bmatrix}
 <d^2T.
\]
\end{lemma}

\begin{proof}
Choose a cut-polynomial certificate $p$ from
\eqref{eq:certificate-equivalence} and write it in the integral-Hessian form
\[
 p(z)=\frac12z^{\top}Hz+\ell^{\top}z+\alpha.
\]
By \eqref{eq:newton-to-hessian}, $H\in\Smat^n(\Z)$, and boundedness below
implies $H\succeq0$.  Since $p$ is integer-valued and
$\min_{\R^n}p>-1$,
\[
 p(\Z^n)\subseteq\Z_+.
\]
Set
\[
 d:=H_{11},
 \qquad
 \sigma_p:=p(0)+p(e_1)\in\Z_+.
\]
The quadratic midpoint identity gives
\begin{equation}\label{eq:hard-midpoint}
 p\!\left(\frac12e_1\right)
 =\frac{p(0)+p(e_1)}2-\frac d8
 =\frac{\sigma_p}{2}-\frac d8.
\end{equation}
Because $\min p>-1$, we obtain
\begin{equation}\label{eq:hard-anchor-bound}
 d<4\sigma_p+8,
 \qquad\text{hence}\qquad
 d\le4\sigma_p+7.
\end{equation}

For a quadratic with Hessian $H$, evaluation at moments with covariance
$\Sigma$ satisfies
\[
 L_p(\bar\eta)=p(\bar x)+\frac12\operatorname{tr}(H\Sigma).
\]
Substituting \eqref{eq:hard-point} and \eqref{eq:hard-midpoint} yields
\begin{equation}\label{eq:hard-anchor-master}
 L_p(\bar\eta)
 =\frac{\sigma_p}{2}
  +\frac{\varepsilon}{2T}
    \bigl(\operatorname{tr}(HG)-dT\bigr).
\end{equation}
The matrices $H$ and $G$ are positive semidefinite, so
$\operatorname{tr}(HG)\ge0$.  If $\sigma_p\ge1$, then
\eqref{eq:hard-anchor-bound} and $\varepsilon=1/44$ give
\[
 L_p(\bar\eta)
 \ge\frac{\sigma_p}{2}-\frac{4\sigma_p+7}{88}
 =\frac{40\sigma_p-7}{88}>0,
\]
contradicting violation.  Hence $\sigma_p=0$, and
nonnegativity on $\Z^n$ gives $p(0)=p(e_1)=0$.
The bound \eqref{eq:hard-anchor-bound} gives $d\le7$.  Positive
semidefiniteness gives $d\ge0$, while $d=0$ would make
\eqref{eq:hard-anchor-master} nonnegative.  Hence
\[
 1\le d\le7.
\]

Let $h\in\Z^m$ contain the last $m$ entries of the first column of $H$.
The Schur complement gives
\[
 H-\frac1d
 \begin{bmatrix}d\\h\end{bmatrix}
 \begin{bmatrix}d\\h\end{bmatrix}^{\top}
 \succeq0.
\]
Taking the trace product with $G\succeq0$ yields
\[
 \operatorname{tr}(HG)
 \ge\frac1d
 \begin{bmatrix}d\\h\end{bmatrix}^{\top}
 G
 \begin{bmatrix}d\\h\end{bmatrix}.
\]
Since $\sigma_p=0$ and $L_p(\bar\eta)<0$,
\eqref{eq:hard-anchor-master} gives
$\operatorname{tr}(HG)<dT$.  Combining the last two inequalities proves
the claim.
\end{proof}

\begin{theorem}[Strong NP-completeness and consequences]
\label{thm:strong-NP-completeness}
The language $\CGViol$ is strongly NP-complete, even for input points in
$\mathcal P_n$. Consequently, exact separation over $\mathcal C_n$ is
strongly NP-hard, and $\CGMem$ is strongly coNP-complete.
\end{theorem}

\begin{proof}
Membership in NP follows from Theorem~\ref{thm:polynomial-certificates}.  We
prove hardness using the construction above.

Suppose first that $D_H\le k$.  Lemma~\ref{lem:scalar-invariance} with $d=1$
gives $h\in\Z^m$ such that, for $v=(1,h)\in\Z^n$,
\[
 v^{\top}Gv=110D_H\le110k<T.
\]
The vector $v$ is primitive because its anchor component is one.  Consider
the primitive split quadratic
\[
 q_v(z):=\frac12(v^{\top}z)(v^{\top}z-1).
\]
Since $q_v$ is integer-valued on $\Z^n$ and
$\min_{z\in\R^n}q_v(z)=-1/8$, it is a rank-one CG cut polynomial by
Lemma~\ref{lem:CG-polynomial-dictionary}.
Since $\bar x=\frac12e_1$ and $v=(1,h)$,
we have $v^\top\bar x=\frac12$. Using
$\bar X=\bar x\bar x^\top+\Sigma$ and~(6.5),
we obtain
\begin{align*}
L_{q_v}(\bar\eta)
&=\frac12(v^\top\bar Xv-v^\top\bar x)\\
&=\frac12\left((v^\top\bar x)^2
  +v^\top\Sigma v-v^\top\bar x\right)=\frac12\left(v^\top\Sigma v-\frac14\right)\\
&=\frac12\left(\frac14-\varepsilon
  +\frac{\varepsilon}{T}v^\top Gv-\frac14\right)\\
&=\frac{\varepsilon}{2T}(v^\top Gv-T)<0,
\end{align*}
where the last inequality follows from
$v^\top Gv<T$.
Thus, $\bar\eta\notin\mathcal C_n$.

Now suppose that $D_H>100k$ and, toward a contradiction, that
$\bar\eta\notin\mathcal C_n$.  Lemma~\ref{lem:arbitrary-cut-extraction}
gives $d\in\{1,\ldots,7\}$ and $h\in\Z^m$ such that
\[
\begin{bmatrix}d\\h\end{bmatrix}^{\top}G\begin{bmatrix}d\\h\end{bmatrix}<d^2T\le7^2\cdot220k=10780k.
\]
However, since $1\le d\le7<11$, Lemma~\ref{lem:scalar-invariance} and \eqref{eq:hard-Gram-distance} imply
\[
\begin{bmatrix}d\\h\end{bmatrix}^{\top}G\begin{bmatrix}d\\h\end{bmatrix}\ge110D_H>11000k,
\]
a contradiction.  Hence $\bar\eta\in\mathcal C_n$ in the NO case.  We have
proved
\[
 D_H\le k\ \Longrightarrow\ \bar\eta\notin\mathcal C_n,
 \qquad
 D_H>100k\ \Longrightarrow\ \bar\eta\in\mathcal C_n,
\]
which has the required orientation for a deterministic Karp reduction to
$\CGViol$.

The entries of $B$ and $\tau$ are bounded by an absolute constant, so $G$
has entries $O(N)$. Since $T=N$ and $n=1+r+10N=\Theta(N)$,
\eqref{eq:hard-point} and the conversion to integral moment coordinates give
polynomial bounds in $n$ for all numerators and denominators of $\bar\eta$.
The construction is polynomial-time and produces points in $\mathcal P_n$,
proving strong NP-completeness.

Finally, $\CGViol$ is the decision version of exact separation.  Hence exact
separation over $\mathcal C_n$ is strongly NP-hard.  Since $\CGMem$ is the
complement of the strongly NP-complete language $\CGViol$, it is strongly
coNP-complete.
\end{proof}

\begin{corollary}[Strong NP-hardness of split-inequality separation]
\label{cor:split-separation-hardness}
When $n$ is part of the input, deciding whether a rational point violates some
split inequality is strongly NP-hard. Consequently, exact separation over
$\mathcal S_n$ is strongly NP-hard.
This resolves the conjecture of Burer and Letchford
\cite[p.~255]{BurerLetchford2014}.
\end{corollary}

\begin{proof}
Use the reduction in the proof of
Theorem~\ref{thm:strong-NP-completeness}.  If $D_H\le k$, the vector
$v=(1,h)$ constructed there is primitive, and the split quadratic
$q_v=\widehat q_{v,0}$ satisfies $L_{q_v}(\bar\eta)<0$.  Hence
$\bar\eta\notin\mathcal S_n$.  If $D_H>100k$, the same proof gives
$\bar\eta\in\mathcal C_n$, and therefore $\bar\eta\in\mathcal S_n$ by
\eqref{eq:universal-chain}.  The reduction has bounded rational data and
produces only points in $\mathcal P_n$, so the claim follows.
\end{proof}

%% file: sections/07_conclusion.tex
\section{Conclusion}
\label{sec:conclusion}

We characterize the facets of $\MIQ_n$ through perfect Delaunay polyhedra
and show that a facet inequality is valid for the rank-one CG closure
$\mathcal C_n$ exactly when its depth, measured after primitive integer-valued
normalization of the defining quadratic, is less than one. Known
classifications and exact depth calculations then establish that one CG
round recovers the hull exactly for $1\le n\le6$, while split inequalities
do so exactly for $1\le n\le5$, with
$\mathcal C_n\subsetneq\mathcal S_n$ for every $n\ge6$. For rational input,
we prove the existence of polynomial-size certificates of nonmembership in
$\mathcal C_n$. When $n$ is part of the input, nonmembership is strongly
NP-complete and membership is strongly coNP-complete. The same reduction
proves strong NP-hardness of quadratic split-inequality separation,
resolving Burer and Letchford's conjecture
\cite[p.~255]{BurerLetchford2014}.

Future work could determine the CG rank of the $ER_7$ facet and develop
criteria based on depth for recovery under iterated closures. Another
direction is to characterize the facets and other maximal proper faces of
mixed-integer quadratic moment hulls, including their inheritance from the
integer moment hull and the effects of nonnegativity or linear coupling
constraints.

%% file: appendices/A_low_dimensional.tex
\appendix
\section{Low-dimensional generators and exact verification}
\label{app:low-dimensional}

We record the exact representatives used in Table~\ref{tab:depths}.  For a quadratic written as
\[
 q(z)=(z-c)^{\top}A(z-c)-r,
\]
the continuous minimum is $-r$, and its depth is therefore $r/\Delta(q)$.

\paragraph{Primitive splits.}
For primitive $v\in\Z^n$ and $s\in\Z$, the quadratic $\widehat q_{v,s}$ of \eqref{eq:split-polynomial} has $\Delta(\widehat q_{v,s})=1$ and continuous minimum $-1/8$.

\paragraph{$2_{21}$ in dimension six.}
The polynomial in \eqref{eq:G6-generator} satisfies
\[
 q_{G_6}(z)
 =\left(z-\frac13\mathbf1\right)^{\top}
   A_{G_6}
  \left(z-\frac13\mathbf1\right)-\frac23,
\]
where
\[
A_{G_6}=\begin{bmatrix}
1&\frac12&0&0&\frac12&-\frac12\\
\frac12&1&\frac12&0&0&-\frac12\\
0&\frac12&1&\frac12&0&-\frac12\\
0&0&\frac12&1&\frac12&-\frac12\\
\frac12&0&0&\frac12&1&-\frac12\\
-\frac12&-\frac12&-\frac12&-\frac12&-\frac12&1
\end{bmatrix}.
\]
Its leading principal minors are
\[
1,\quad\frac34,\quad\frac12,\quad\frac5{16},\quad
\frac18,\quad\frac3{64},
\]
so $A_{G_6}\succ0$.  The nonconstant Newton coefficients have greatest common divisor one, and hence $\delta(2_{21})=2/3$.

\paragraph{$3_{21}$ in dimension seven.}
Let
\[
 A_7=\begin{bmatrix}
 4&3&3&3&3&5&4\\
 3&4&3&3&3&5&4\\
 3&3&4&3&3&5&4\\
 3&3&3&4&3&5&4\\
 3&3&3&3&4&5&4\\
 5&5&5&5&5&8&6\\
 4&4&4&4&4&6&6
 \end{bmatrix},
 \qquad
 c_7=\frac12e_7.
\]
Then
\[
 A_7c_7=(2,2,2,2,2,3,3)^{\top},
 \qquad
 c_7^{\top}A_7c_7=\frac32,
\]
and
\[
 q_{G_7}(z):=(z-c_7)^{\top}A_7(z-c_7)-\frac32
\]
has zero constant term.  Its Newton coefficients are
\[
 a=(0,0,0,0,0,2,0),
 \qquad
 d=(8,8,8,8,8,16,12),
 \qquad
 \gamma_{ij}=2(A_7)_{ij}.
\]
Their greatest common divisor is two, so
\[
 \widehat q_{G_7}=\frac12q_{G_7},
 \qquad
 \delta(3_{21})=\frac34.
\]
The leading principal minors of $A_7$ are
\[
4,\quad7,\quad10,\quad13,\quad16,\quad3,\quad2,
\]
so $A_7\succ0$.

\paragraph{$ER_7$.}
Let
\[
 A_{ER}=\begin{bmatrix}
 8&6&6&6&6&6&9\\
 6&8&6&6&6&6&9\\
 6&6&8&6&6&6&9\\
 6&6&6&8&6&6&9\\
 6&6&6&6&8&6&9\\
 6&6&6&6&6&8&9\\
 9&9&9&9&9&9&13
 \end{bmatrix},
 \qquad
 c_{ER}=\frac1{16}(5,5,5,5,5,5,-14)^{\top}.
\]
Then
\[
 A_{ER}c_{ER}=
 \left(4,4,4,4,4,4,\frac{11}{2}\right)^{\top},
 \qquad
 c_{ER}^{\top}A_{ER}c_{ER}=\frac{43}{16}.
\]
The polynomial
\[
 q_{ER}(z):=(z-c_{ER})^{\top}A_{ER}(z-c_{ER})-\frac{43}{16}
\]
has Newton coefficients
\[
 a=(0,0,0,0,0,0,2),
 \qquad
 d=(16,16,16,16,16,16,26),
 \qquad
 \gamma_{ij}\in\{12,18\}.
\]
Their greatest common divisor is two, so
\[
 \widehat q_{ER}=\frac12q_{ER},
 \qquad
 \delta(ER_7)=\frac{43}{32}.
\]
The leading principal minors of $A_{ER}$ are
\[
8,\quad28,\quad80,\quad208,\quad512,\quad1216,\quad256,
\]
so $A_{ER}\succ0$.

%% file: appendices/B_encoding_bounds.tex
\section{Encoding bounds for separation certificates}
\label{app:encoding-bounds}

This appendix supplies the quantitative bit-length estimates used in the proof of Theorem~\ref{thm:polynomial-certificates}.  All input rationals are encoded by binary numerators and positive binary denominators.  

\subsection{Indefinite moment matrices}
\label{subsec:encoding-indefinite}

Let $\bar W\in\Q^{s\times s}$ be symmetric and indefinite, where $s=n+1$, and let $L$ denote its binary input length.  Clear denominators and write
\[
 S=D\bar W\in\Z^{s\times s},
 \qquad D>0.
\]
Taking $D$ to be the product of the input denominators gives
$\log D=O(L)$.  Put
\[
 M:=\max\{1,|S_{ij}|:1\le i,j\le s\}.
\]
Then $\log M$ is polynomial in $L$.

Every eigenvalue $\mu$ of $S$ satisfies $|\mu|\le sM$.  Indeed, if $Sx=\mu x$ and $|x_i|=\|x\|_\infty$, then
\[
 |\mu|\,|x_i|
 \le\sum_{j=1}^s|S_{ij}|\,|x_j|
 \le sM|x_i|.
\]
Let $\rho=\rank(S)$ and let $\mu_1,\ldots,\mu_\rho$ be the nonzero eigenvalues.  Up to sign, their product is the coefficient of $t^{s-\rho}$ in the integer polynomial $\det(tI-S)$, and hence has absolute value at least one. Therefore every nonzero eigenvalue has magnitude at least
\[
 \Delta:=(sM)^{-(s-1)}.
\]
Since $S$ is indefinite, it has a unit eigenvector $u$ whose eigenvalue $\lambda$ satisfies $\lambda\le-\Delta$.

Choose the integer
\[
 K:=8s^2M\Delta^{-1}=8s^2M(sM)^{s-1}
\]
and round $Ku$ coordinatewise to $v\in\Z^s$.  With $y:=v/K$,
\[
 \|y-u\|_2\le\frac{\sqrt s}{2K}.
\]
Because $\|S\|_2\le sM$, $\|u\|_2=1$, and
$\|y\|_2\le1+\sqrt s/(2K)$,
\[
 \begin{aligned}
 |y^{\top}Sy-u^{\top}Su|
 &\le\|S\|_2\,\|y-u\|_2(\|y\|_2+\|u\|_2)\\
 &\le\frac{3s^{3/2}M}{2K}
 \le\frac{3\Delta}{16}
 <\frac{\Delta}{2}.
 \end{aligned}
\]
Thus
\[
 y^{\top}Sy< -\frac{\Delta}{2}<0,
 \qquad
 v^{\top}Sv=K^2y^{\top}Sy<0.
\]
Since $D>0$, this is equivalent to $v^{\top}\bar Wv<0$.

Finally,
\[
 \log\Delta^{-1}=(s-1)(\log s+\log M)=\operatorname{poly}(L),
\]
so $\log K=\operatorname{poly}(L)$.  Coordinatewise rounding gives $|v_i|\le K+1$, and hence the binary length of $v$ is polynomial in $L$. The Newton coefficients of $(v_0+v'^{\top}z)^2$ are sums and pairwise products of entries of $v$, so their total binary length is polynomial as well.

\subsection{Positive-definite covariance}
\label{subsec:encoding-positive-definite}

Let $\bar W\succeq0$ be rational, let $\Sigma=\bar X-\bar x\bar x^{\top}\succ0$, and let
\[
 p(z)=z^{\top}Az+2b^{\top}z+\alpha
\]
be any cut-polynomial certificate for $\bar\eta$.  As shown in \eqref{eq:positive-covariance-master}--\eqref{eq:positive-covariance-trace}, if $\min p=-\delta$, then $0<\delta<1$ and
\begin{equation}
\label{eq:appendix-trace-bound}
 \ip{A}{\Sigma}<1.
\end{equation}

Write
\[
 \Sigma=\frac{1}{D_\Sigma}S,
 \qquad
 D_\Sigma\in\Z_+,
 \quad
 S\in\Z^{n\times n}\text{ positive definite},
\]
and let
\[
 M_\Sigma:=\max\{1,|S_{ij}|:1\le i,j\le n\}.
\]
The numbers $D_\Sigma$ and $M_\Sigma$ have polynomial binary length in the input size.  Since $\det S$ is a positive integer and $\lambda_{\max}(S)\le nM_\Sigma$,
\[
 \lambda_{\min}(\Sigma)
 =\frac{\lambda_{\min}(S)}{D_\Sigma}
 \ge\frac{1}{D_\Sigma(nM_\Sigma)^{n-1}}.
\]
Define
\[
 B_\Sigma:=D_\Sigma(nM_\Sigma)^{n-1}.
\]
The binary length of $B_\Sigma$ is polynomial in the input length.  From \eqref{eq:appendix-trace-bound},
\begin{equation}
\label{eq:appendix-A-bound}
 \|A\|_2\le\operatorname{tr}(A)<B_\Sigma.
\end{equation}

Choose $y$ with $Ay=-b$.  The master identity gives $(\bar x-y)^{\top}A(\bar x-y)<1$, and therefore
\[
 \|A\bar x+b\|_2^2
 =(\bar x-y)^{\top}A^2(\bar x-y)
 \le\|A\|_2(\bar x-y)^{\top}A(\bar x-y)
 <B_\Sigma.
\]
Put
\[
 X_\infty:=\max\{1,|\bar x_i|:1\le i\le n\}.
\]
This rational number has polynomial binary length, and $\|\bar x\|_2\le\|\bar x\|_1\le nX_\infty$.  Since $B_\Sigma\ge1$, we have $\sqrt{B_\Sigma}\le B_\Sigma$, and hence
\begin{equation}
\label{eq:appendix-b-bound}
 \|b\|_\infty\le\|b\|_2
 \le B_\Sigma(1+nX_\infty).
\end{equation}
Since $p$ is integer-valued and its continuous minimum is greater than $-1$, $p(0)=\alpha\in\Z_+$.  Moreover,
$L_p(\bar\eta)<0$ and $\ip{A}{\Sigma}\ge0$ imply $p(\bar x)<0$, whence
\begin{equation}
\label{eq:appendix-alpha-bound}
 \begin{aligned}
 0\le\alpha
 &< -\bar x^{\top}A\bar x-2b^{\top}\bar x\\
 &\le B_\Sigma n^2X_\infty^2
    +2B_\Sigma(1+nX_\infty)nX_\infty.
 \end{aligned}
\end{equation}
All quantities on the right-hand sides of \eqref{eq:appendix-A-bound}--\eqref{eq:appendix-alpha-bound} have polynomial binary length. By \eqref{eq:newton-to-matrix}, the entries of $A$ and $b$ have denominators at most two and four, respectively, while $\alpha\in\Z$. The magnitude bounds therefore give polynomial bounds on their binary encoding lengths. The Newton coefficients are recovered by
\[
 d_i=2A_{ii},
 \qquad
 \gamma_{ij}=2A_{ij},
 \qquad
 a_i=2b_i+A_{ii},
\]
so they also have polynomial total encoding length.

\subsection{Singular covariance}
\label{subsec:encoding-singular}

Let $\Sigma\succeq0$ be rational and singular.  Clear denominators and write
$S=D_\Sigma\Sigma\in\Z^{n\times n}$.  Put $r:=\rank(S)$ and $k:=n-r$.
When $r>0$, select $r$ linearly independent rows of $S$ and collect them in
$A_0\in\Z^{r\times n}$.  The full-row-rank extension of the Kannan--Bachem
Hermite-normal-form algorithm
\cite[Section~2, Theorem~2 and p.~504]{KannanBachem} computes
\begin{equation}\label{eq:appendix-HNF}
 A_0Q=\begin{bmatrix}H&0\end{bmatrix},
 \qquad
 Q\in\GL_n(\Z),
 \quad
 H\in\Z^{r\times r}\text{ nonsingular},
\end{equation}
with polynomial running time and polynomially bounded output bit length.  If
$r=0$, let $A_0$ be the empty $0\times n$ matrix and take $Q=I_n$.  In either
case, let $Z$ be the last $k$ columns of $Q$.  The selected rows span the row
space of $S$, so
$\ker_{\R}(A_0)=\ker_{\R}(S)=\ker_{\R}(\Sigma)$.  Moreover, for
$x=Qy\in\Z^n$, equation~\eqref{eq:appendix-HNF} gives
\[
 Sx=0
 \quad\Longleftrightarrow\quad
 A_0x=0
 \quad\Longleftrightarrow\quad
 Hy_{1:r}=0
 \quad\Longleftrightarrow\quad
 y_{1:r}=0.
\]
Thus the columns of $Z$ form a lattice basis of
$\ker(\Sigma)\cap\Z^n$.  Set
\[
 R:=Z^{\top},
 \qquad
 U:=Q^{-\top}.
\]
Because $\Sigma$ is symmetric,
$\range(\Sigma)=\ker(\Sigma)^\perp$, and the identity
$Q^{\top}Q^{-\top}=I_n$ yields
\begin{equation}\label{eq:appendix-RU}
 \ker_{\R}R=\range(\Sigma),
 \qquad
 RU=\begin{bmatrix}0&I_k\end{bmatrix}.
\end{equation}
The matrix $R$ has polynomial binary length.  Since $Q$ is unimodular, the
cofactor formula gives the same conclusion for $U=Q^{-\top}$, while
$U^{-1}=Q^{\top}$ already has polynomial binary length.

Suppose first that $R\bar x\notin\Z^k$.  Choose a row $v^{\top}$ of $R$
with $v^{\top}\bar x\notin\Z$ and set
$s=\lfloor v^{\top}\bar x\rfloor$.  The vector $v$ is primitive because it
is the transpose of a column of the unimodular matrix $Q$, and
$v\in\ker(\Sigma)$.  With
$\theta=v^{\top}\bar x-s\in(0,1)$,
\[
 L_{\widehat q_{v,s}}(\bar\eta)
 =\frac12\left[v^{\top}\Sigma v+\theta(\theta-1)\right]
 =\frac12\theta(\theta-1)<0.
\]
The entries of $v$ and the integer $s$ have polynomial binary length, so the
split polynomial is a polynomial-size certificate.

Suppose now that $R\bar x\in\Z^k$.  Define
\begin{equation}\label{eq:appendix-z0}
 z^0:=U
 \begin{bmatrix}0\\R\bar x\end{bmatrix}\in\Z^n.
\end{equation}
Then $R(z^0-\bar x)=0$, and hence
$z^0-\bar x\in\range(\Sigma)$ by \eqref{eq:appendix-RU}.  Apply the
integer-affine transformation
\[
 z=z^0+Uy.
\]
The transformed mean and covariance are
\begin{equation}\label{eq:appendix-transformed-moments}
 \bar x'=U^{-1}(\bar x-z^0),
 \qquad
 \Sigma'=U^{-1}\Sigma U^{-\top}=Q^{\top}\Sigma Q.
\end{equation}
Equation~\eqref{eq:appendix-RU} shows that the last $k$ coordinates of
$\bar x'$ vanish.  Partition $Q=[Q_1\ Z]$.  Since $\Sigma Z=0$,
\[
 \Sigma'
 =\begin{bmatrix}
 Q_1^{\top}\Sigma Q_1&0\\0&0
 \end{bmatrix}
 =:\begin{bmatrix}\Sigma_r&0\\0&0\end{bmatrix}.
\]
Congruence preserves rank, so $\Sigma_r$ is an $r\times r$ positive-definite
matrix.  The transformed moment point is therefore the exact zero extension
of an $r$-dimensional rational point $\bar\eta_r$ with positive-definite
covariance.  If $r=0$, then $\bar x'=0$ and the transformed point is the
integer generator $W_0$, contradicting nonmembership.  Hence $r\ge1$.

All entries of $\bar\eta_r$ have polynomial binary length by
\eqref{eq:appendix-HNF}--\eqref{eq:appendix-transformed-moments}.  The
positive-definite case supplies a certificate $p_r$ of polynomial size for
$\bar\eta_r$.  Extend $p_r$ so that it is independent of the last $k$
transformed variables and pull it back by
\begin{equation}\label{eq:appendix-pullback}
 p(z):=p_r\!\left(
 \left[U^{-1}(z-z^0)\right]_{1:r}
 \right).
\end{equation}
Integer-affine invariance shows that $p$ is integer-valued, has the same
continuous minimum as $p_r$, and satisfies
$L_p(\bar\eta)=L_{p_r}(\bar\eta_r)<0$.  The entries of $U^{-1}$ and $z^0$
have polynomial binary length, and expanding the quadratic in
\eqref{eq:appendix-pullback} requires only polynomially many arithmetic
operations on polynomial-size numbers.  The resulting Newton coefficients
therefore have polynomial total encoding length.